\documentclass[12pt, notitlepage]{amsart}
\usepackage{latexsym, amsfonts, amsmath, amssymb, amsthm, cite}

\usepackage{stmaryrd}

\newtheorem{theorem}{Theorem}[section]

\newtheorem{conjecture}{Conjecture}

\usepackage{graphicx}
\usepackage{mathrsfs}
\numberwithin{equation}{section}

\newcommand{\E}{\mathbb E}
\begin{document}
 \title{The distribution of values of zeta and L-functions}
%\titlemark{Zeta and L-functions}

\author{Kannan Soundararajan} 
 \address{Department of Mathematics, Stanford University} 
  \email{ksound@stanford.edu}

%\dedication{Dedicated to ...}

%------
% Insert your abstract.
%------
%\begin{abstract}
%We survey recent progress on understanding the distribution of values of zeta and $L$-functions.  In particular, we discuss the problem of moments of $|\zeta(\tfrac 12+it)|$ and moments of central $L$-values in families, where the last twenty five years have seen a conjectural understanding of the asymptotics of these moments, together with progress in obtaining good upper and lower bounds in many situations.   

%\end{abstract}

\maketitle

%------
% INSERT THE BODY OF THE PAPER HERE (except
% acknowledgments, funding info and bibliography)
%

\noindent This article concerns the distribution of values of the Riemann zeta-function, and related $L$-functions.   
We begin with a brief discussion of $L$-values at the edge of the critical strip, which give information on arithmetic invariants such as class 
numbers.  The remainder of the article is concerned with the value distribution of $\zeta(\tfrac 12+it)$ and the distribution of central values in 
families of $L$-functions.   The typical behavior of $\zeta(\tfrac 12+it)$ is described by a fundamental theorem of Selberg (discussed in \S 2) which asserts that 
$\log \zeta(\tfrac 12+it)$ is distributed like a complex Gaussian with prescribed mean and variance.   Analogues of Selberg's theorem for central values in families of $L$-functions were conjectured by Keating and Snaith, and we motivate these conjectures and the progress towards them in \S 3.    Section 4 begins our treatment of the problem of understanding the moments 
of $|\zeta(\tfrac 12+it)|$ and analogous questions for central $L$-values.    While this is a classical topic, going back to work of Hardy and Littlewood, 
it is only in the last twenty five years that even a good conjectural understanding of the problem has emerged.  The Keating--Snaith conjectures for the asymptotics 
of moments were first developed by pursuing an analogy between values of the zeta function and the values of the characteristic polynomial of large random matrices.   
These conjectures are described in \S 5, which also shows how the problem of understanding moments is tied up with understanding the \emph{large deviations} range in Selberg's theorem.  Progress towards the moment conjectures (see \S 6) has been of three types: (i) understanding asymptotics for small moments in a number of examples, (ii) obtaining lower bounds of the correct order of magnitude (which are known in many cases), and (iii) obtaining in great generality upper bounds of the correct order of magnitude assuming the Generalized Riemann Hypothesis.   In \S 7 we discuss what is known about the maximal size of $|\zeta(\tfrac 12+it)|$ and 
central $L$-values, and speculate on what the truth might be.   Finally, in \S 8 we consider briefly an intriguing problem of Fyodorov--Hiary--Keating on understanding the ``local maximum'' of $|\zeta(\tfrac 12+it)|$ for $t$ in intervals of length $1$, which is  closely connected to problems in \emph{branching Brownian motion} and \emph{Gaussian multiplicative chaos}.

\section{Values at the edge of the critical strip}

\noindent It was already observed by Gauss and Dirichlet that certain special values of $L$-functions encode interesting arithmetic information.  
Recall that a \emph{discriminant} is an  integer $d \equiv 0, 1 \pmod 4$, and $d$ is called a \emph{fundamental discriminant} if $d/m^2$ is not a discriminant for any divisor $m^2$ of $d$ larger than $1$.    Fundamental discriminants are in one-to-one correspondence with discriminants of quadratic fields ${\mathbb Q}(\sqrt{d})$.  Associated to a fundamental discriminant $d$ is the Kronecker--Legendre symbol $\chi_d(n) = (\frac{d}{n})$, which is a primitive Dirichlet character $\pmod {|d|}$.  For example, if $p$ is an odd prime then either $p$ or $-p$ is a fundamental discriminant (depending on whether $p$ is $1$ or $3\pmod 4$), and in either case the associated quadratic character is the familiar Legendre symbol $(\bmod \, p)$.    Associated to the primitive character $\chi_d$ is the Dirichlet $L$-function 
$$ 
L(s,\chi_d) = \sum_{n=1}^{\infty} \frac{\chi_d(n)}{n^s} = \prod_p \Big( 1- \frac{\chi_d(p)}{p^s}\Big)^{-1}. 
$$ 
Although $d=1$ is permitted in our definition of fundamental discriminants (and corresponds to the Riemann zeta-function), it is an anomalous case and we 
shall mainly be interested in fundamental discriminants $d\neq 1$.  Like the Riemann zeta-function, the Dirichlet $L$-function $L(s,\chi_d)$ converges absolutely for Re$(s) >1$, extends analytically to the entire complex plane (unlike $\zeta(s)$, there is no pole at $s=1$ here), and satisfies a functional equation connecting values at $s$ to values at $1-s$.  The non-trivial zeros of $L(s,\chi_d)$ lie in the \emph{critical strip} $0< \text{Re}(s) <1$, with the Generalized Riemann Hypothesis (GRH) predicting that they lie on the critical line Re$(s)= \frac 12$.   For background on Dirichlet $L$-functions see Davenport \cite{Dav}, and for a general 
comprehensive treatment of analytic number theory (including information on many other families of $L$-functions that will be considered here) see Iwaniec and Kowalski 
\cite{IK}.

In this family of quadratic Dirichlet $L$-functions, the values $L(1,\chi_d)$ (lying at the edge of the critical strip) are of great arithmetical interest.  A key step in Dirichlet's proof that there are infinitely many primes in arithmetic progressions involves showing that $L(1,\chi_d) \neq 0$.   Dirichlet established this by finding a beautiful connection between $L(1,\chi_d)$ and the group of equivalence classes of binary quadratic forms of discriminant $d$ which had earlier been studied by Gauss.   For example, if $d$ is a negative fundamental discriminant, then Dirichlet's class number formula states that 
$$ 
L(1,\chi_d) = \frac{2\pi}{w}  \frac{h(d)}{\sqrt{|d|}}, 
$$ 
where $h(d)$ is  a positive integer, namely the class number of the imaginary quadratic field ${\mathbb Q}(\sqrt{d})$, and $w$ counts the number of roots of unity in ${\mathbb Q}(\sqrt{d})$ (so that $w=2$ for $d<-4$, and $w=4$ for $d=-4$, and $w=6$ for $d=-3$).   
The special case $d=-4$ of the Dirichlet class number formula is widely familiar as the Madhava--Leibniz--Gregory series $1-1/3+1/5-1/7+\ldots =\pi/4$.   Another classical connection to these special $L$-values arises in the Gauss--Legendre three squares theorem.   If $n$ is a square-free integer with $n\equiv 3 \pmod 8$, then the number of ways of writing $n$ as a sum of three squares, $r(n)$, equals $24 h(-n)$; a result known to Gauss, together with variants when $n\equiv 1, 2 \pmod 4$.  

These connections motivate the study of the distribution of the values $L(1,\chi_d)$.  Here are some natural questions that arise.  If fundamental discriminants $d$ are chosen uniformly with $|d|\le X$, (i) what is the statistical distribution of the values $L(1,\chi_d)$, and (ii) what are the largest and smallest possible values of $L(1,\chi_d)$?   As we shall see, the problem of the statistical distribution of $L(1,\chi_d)$ can be understood quite precisely, but there are still large gaps in our understanding of the extreme values. 

Let us begin with the simpler situation of $L(2,\chi_d)$ where both the Dirichlet series and Euler product in the definition of $L(s,\chi_d)$ converge absolutely.  If the values $\chi_d(p)$ are known for all primes $p\le z$ then 
$$ 
\Big| L(2,\chi_d) - \prod_{p\le z} \Big(1-\frac{\chi_d(p)}{p^2} \Big)^{-1} \Big| \le \sum_{n>z} \frac{1}{n^2} = O\Big( \frac1z\Big).
$$ 
The value $\chi_d(p) = (\frac{d}{p})$ is determined by $d\pmod p$ for odd primes $p$, and for $p=2$ the value of $\chi_d(p)$ is determined by $d\pmod 8$.  Thus, by the Chinese Remainder Theorem, the values of $\chi_d(p)$ for $p\le z$ are determined by $d$ modulo $4\prod_{p\le z} p$.  
One way to view this is as a kind of \emph{almost periodicity}:  if two fundamental discriminants $d_1$ and $d_2$ are congruent modulo $4\prod_{p\le z} p$ then $L(2,\chi_{d_1}) =L(2,\chi_{d_2}) + O(1/z)$.

If $p$ is an odd prime and $X$ is large, then a little calculation shows that a proportion $\tfrac{1}{p+1}$ of the fundamental discriminants $d$ with $|d|\le X$ are multiples of $p$ (this is essentially the proportion of square-free integers that are multiples of $p$) and $\chi_d(p)=0$ here.  The remaining proportion $\frac{p}{p+1}$ of fundamental discriminants are evenly split among the possible values $\chi_d(p) = 1$ or $-1$.  Pleasantly, it turns out that for $p=2$ also a proportion $\tfrac 13$ of the fundamental discriminants $|d|\le X$ satisfy each of the three cases $\chi_d(2) = 0$, $1$ or $-1$.   Moreover the Chinese Remainder Theorem tells us that for different primes $p$, the values $\chi_d(p)$ are distributed ``independently'' of each other, at least if we restrict to primes $p\le z$ with $\prod_{p\le z} p$ being small in comparison with $X$.  This motivates us to define for prime numbers $p$, independent random variables ${\mathbb X}(p)$ taking the values $0$ with probability $1/(p+1)$ and the values $\pm 1$ with probability $p/(2(p+1))$.   Then the distribution of $\prod_{p\le z}(1-\chi_d(p)/p^2)^{-1}$ is the same as the distribution of the random Euler product $\prod_{p\le z} (1- {\mathbb X}(p)/p^2)^{-1}$.  Letting $z\to \infty$, we have described the distribution of $L(2,\chi_d)$ as being precisely the distribution of $\prod_p (1- {\mathbb X}(p)/p^2)^{-1}$.  

The story for extreme values is also clear:  
$$ 
\frac{\zeta(4)}{\zeta(2)} = \prod_p \Big(1 + \frac 1{p^2}\Big)^{-1} \le \prod_p \Big(1- \frac{\chi_d(p)}{p^2}\Big)^{-1} = L(2,\chi_d) \le 
\prod_p \Big(1 -\frac{1}{p^2}\Big)^{-1} = \zeta(2). 
$$ 
Moreover we may find values $L(2,\chi_d)$ arbitrarily close to $\zeta(4)/\zeta(2)$ by choosing $d$ with $\chi_d(p) =-1$ for all primes $p\le z$, and we may find values arbitrarily close to $\zeta(2)$ by choosing $d$ with $\chi_d(p) =1$ for all primes $p\le z$.

Let us now turn to the distribution of $L(1,\chi_d)$ where there is a similar story but with some added complications since the series and product defining $L(s,\chi_d)$ are no longer absolutely convergent.   For example, one can show that if $z\le (\log X)^{10}$ then $L(1,\chi_d) = \prod_{p\le z} (1-\chi_d(p)/p)^{-1} + O(1/z^{\frac 14})$ for all but $O(X/z^{\frac 14})$ of the fundamental discriminants $|d|\le X$.  This again may be viewed as a kind of almost periodicity:  allowing $z$ to tend slowly to infinity with $X$, for almost all pairs of discriminants $d_1$ and $d_2$ with $d_1 \equiv d_2 \pmod{4\prod_{p\le z}p}$ one has $L(1,\chi_{d_1})\approx L(1,\chi_{d_2})$.

For primes $p$, let ${\mathbb X}(p)$ denote the random variables described earlier, and extend ${\mathbb X}$ to all integers using (complete) multiplicativity; thus, if $n= p_1^{e_1}\cdots p_k^{e_k}$ then ${\mathbb X}(n) = {\mathbb X}(p_1)^{e_1} \cdots {\mathbb X}(p_k)^{e_k}$.   This is an example of a \emph{random multiplicative function}, and we may correspondingly consider the random $L$-function
\begin{equation} 
\label{1.1} 
L(s,{\mathbb X}) =\sum_{n=1}^{\infty} \frac{{\mathbb X}(n)}{n^s} =  \prod_p \Big( 1- \frac{{\mathbb X}(p)}{p^s}\Big)^{-1}. 
\end{equation} 
Both the series and product above converge almost surely provided Re$(s) > \tfrac 12$; this follows essentially from the fact that the variance of 
$\sum_p {\mathbb X}(p)/p^{s}$ is $\sum_{p} \frac{p}{p+1} \frac{1}{p^{2\text{Re}(s)}}$, which is a convergent sum when Re$(s)>\tfrac 12$.  In 
particular, the random Euler product $L(1,{\mathbb X})$ converges almost surely,   and the values $L(1,\chi_d)$ are distributed like  
$L(1,{\mathbb X})$.    We may see this by first approximating most $L(1,\chi_d)$ by $\prod_{p\le z}(1-\chi_d(p)/p)^{-1}$, noting that this truncated 
Euler product is distributed exactly like $\prod_{p\le z}(1-{\mathbb X}(p)/p)^{-1}$, and finally letting $z\to \infty$.

Let us state the result discussed above more precisely.  Given any $\tau >0$, the proportion of fundamental discriminants $|d|\le X$ with 
$L(1,\chi_d) \ge e^{\gamma} \tau$ tends as $X\to \infty$ to $\text{Prob}(L(1,{\mathbb X}) > e^{\gamma} \tau)$.  Here $\gamma$ is Euler's constant, 
and we have normalized in this fashion in view of Mertens's theorem $\prod_{p\le z} (1-1/p)^{-1} \sim e^{\gamma}\log z$.  If $\tau$ is large, and we 
seek values of $L(1,\chi_d)$ larger than $e^{\gamma}\tau$, the most likely way in which such large values arise is when $\chi_d(p) =1$ for all primes 
up to about $e^{\tau}$.  Similarly, the proportion of fundamental discriminants $|d|\le X$ with $L(1,\chi_d) <  \zeta(2)/(e^{\gamma}\tau)$ tends as 
$X\to \infty$ to $\text{Prob}(L(1,{\mathbb X}) <\zeta(2)/( e^{\gamma} \tau))$.   The normalization here is made in view of $\prod_{p\le z} (1+1/p)^{-1} 
\sim \zeta(2)/(e^{\gamma}\log z)$.   The distribution of $L(1,{\mathbb X})$ is continuous --- it is more natural to think of the distribution of $\log L(1,
{\mathbb X})$ which is smooth --- and its tails $\text{Prob}(L(1,{\mathbb X}) > e^{\gamma} \tau)$ or $\text{Prob}(L(1,{\mathbb X}) < \zeta(2)/
(e^{\gamma} \tau))$ decay double exponentially, behaving like $\exp(-(1+o(1))e^{\tau -C_1}/\tau)$ for a suitable constant $C_1$ (see \cite{GS1}).    
With high likelihood one has $1/10 \le L(1,{\mathbb X}) \le 10$, although there is a small positive probability of finding arbitrarily large or arbitrarily 
small values.

The qualitative results mentioned above were obtained by Chowla and Erd{\H o}s \cite{ChoErd}, and with some uniformity in $\tau$ by 
Elliott \cite{Ell1}.  The question of uniformity in $\tau$ is studied in more detail by Montgomery and Vaughan \cite{MV}, and Granville and Soundararajan \cite{GS1}, with the aim of understanding the extreme values of $L(1,\chi_d)$.  By ``uniformity in $\tau$'', we mean the problem of allowing $\tau$ to depend on $X$ while still guaranteeing that the proportion of $|d|\le X$ with $L(1,\chi_d) > e^{\gamma} \tau$ is comparable to the tail probability that $L(1,{\mathbb X}) > e^{\gamma} \tau$ (and similarly for small values of $L(1,\chi_d)$).   In view of the double exponential decay of the tails of the distribution of $L(1, {\mathbb X})$ mentioned above, the largest viable range for uniformity in $\tau$ is $\tau \le  \tau_{\text{max}} +\epsilon$, with $\tau_{\max} =\log\log X+ \log \log \log X + C_1$ and any fixed $\epsilon >0$  --- at this point one has $\text{Prob}(L(1,{\mathbb X}) > e^{\gamma} \tau_{\max}) < 1/X$.  The results in \cite{GS1} show excellent agreement between the distribution of $L(1,\chi_d)$ and the probabilistic model $L(1,{\mathbb X})$ in almost the entire viable range.  These results suggest the following conjectures on the extreme values of $L(1,\chi_d)$: 
\begin{equation} 
\label{1.2} 
\max_{|d|\le X} L(1,\chi_d) = e^{\gamma}  (\tau_{\max} +o(1)), \text{  and   } \min_{|d|\le X} L(1,\chi_d) = \zeta(2)/(e^{\gamma} (\tau_{\max}+ o(1))). 
\end{equation} 
In \cite{GS1} it is shown that there are values of $L(1,\chi_d)$ nearly as large as the conjecture in \eqref{1.2} (for example, assuming the truth of GRH one can find values as large as $e^{\gamma}(\tau_{\max} - C)$ for some constant $C$) and values almost as small as in \eqref{1.2}.  However, as we shall discuss next, there are large gaps in our understanding of why the extreme values cannot be much larger or smaller.

How large can $z$ be such that for some fundamental discriminant $|d|\le X$ one has $\chi_d(p) =1$ for all primes $p\le z$?  This problem is intimately related to finding large values of $L(1,\chi_d)$.  Correspondingly, the problem of finding small values of $L(1,\chi_d)$ may be thought of as wanting $\chi_d(p)=-1$ for all primes $p\le z$.   We noted already that the values of $\chi_d(p)$ for $p\le z$ may be determined by knowing $d\pmod{4\prod_{p\le z} p}$.  The prime number theorem gives $\prod_{p\le z}p = e^{z(1+o(1)}$, so that with $z = \frac 12\log X$ (say) we can find $|d|\le X$ with any given signs $\chi_d(p)$ for $p\le z$ --- for example we may make them all $1$, or all $-1$.  If we think of the probabilistic model ${\mathbb X}$ which treats $\chi_d(p)$ as essentially being a ``coin toss" we may expect that the primes up to about  $z =\log X \log \log X$ (there are about $\log X$ primes below this $z$) may take any prescribed signs.  This dovetails nicely with the conjectured size of extreme values in \eqref{1.2}, since (in the case of large values) $\prod_{p\le z} (1-1/p)^{-1} \sim e^{\gamma} \log z \approx e^{\gamma}(\log \log X +\log \log \log X)$.  For primes $p$ larger than about $\log X \log \log X$, we expect randomness to kick in, and to find an equal number of positive and negative values of 
$\chi_d(p)$.  

Our current knowledge is very far from these probabilistic considerations.  Given a prime $\ell$, Vinogradov conjectured that the least quadratic non-residue $(\bmod \,\ell)$ lies below $C(\epsilon) \ell^{\epsilon}$ for some constant $C(\epsilon)$.   That is, there must be a prime $p \le C(\epsilon) \ell^{\epsilon}$ with $(\frac{p}{\ell}) = -1$, which is a weak version of the prediction from the random model that there exists such $p$ with $p\le C \log \ell \log \log \ell$ for some constant $C$.  Toward Vinogradov's conjecture, we know, as a consequence of the Burgess bounds for character sums, that the least quadratic non-residue lies below $\ell^{1/(4\sqrt{e})+o(1)}$ (see \cite{Bur}), and no improvement over this exponent has been made in more than fifty years.  In terms of $L(1,\chi_d)$, the work towards Vinogradov's conjecture may be used to show that (see \cite{GS2, Ste})
$$ 
L(1,\chi_d) < \Big( \frac{1}{4} \Big(2- \frac{2}{\sqrt{e}}\Big) +o(1) \Big) \log |d|.   
$$ 
This is far from the conjecture in \eqref{1.2}, and even an improvement in the constant above would be significant and lead to an improvement on the bound for the least quadratic non-residue (see also \cite{BG, Tao, GS3} for related work).  

Even less is known about the problem of bounding the least prime $p$ such that $p$ is a quadratic residue $\pmod \ell$.  To give a sense of the interest of this problem, we note that if $\ell \equiv 3\pmod 4$ is a prime, then the imaginary quadratic field ${\mathbb Q}(\sqrt{-\ell})$ has class number $1$ if and only if $(\frac{p}{\ell})=-1$ for all $p< (1+\ell)/4$.  For such a prime $\ell$, the polynomial $n^2+n + (1+\ell)/4$ takes prime values for $0\le n< (\ell-3)/4$.  Euler's famous polynomial $n^2+n+41$ is the largest example of this phenomenon, corresponding to the prime $\ell=163$ for which the first $12$ primes (the primes below 41) are all quadratic non-residues. Toward this problem, we know that the least prime quadratic residue $(\bmod \, \ell)$ lies below $C(\epsilon)\ell^{\frac 14+\epsilon}$ for any $\epsilon >0$ (see \cite{HB6}), but with a constant $C(\epsilon)$ that is \emph{ineffective} (meaning the proof only shows the existence of $C(\epsilon)$, but without any way to compute it, even in principle).   This is related to Siegel's ineffective lower bound (see \cite{Dav}): for any $\epsilon >0$ there exists $C(\epsilon)>0$ with  
$$ 
L(1,\chi_d) > C(\epsilon) |d|^{-\epsilon}.
$$ 
Thus our knowledge of small values of $L(1,\chi_d)$ is even further from the conjecture in \eqref{1.2}.

If we assume the truth of GRH, then much better results are known.   On GRH, the least quadratic non-residue $\pmod \ell$ can be shown to be $< (\log \ell)^2$, and the least prime quadratic residue also lies below $(1+o(1)) (\log \ell)^2$ (see \cite{LaLiSo}).   Moreover, for any fundamental discriminant $d$ one has 
\begin{equation} 
\label{1.3} 
L(1,\chi_d) \sim \prod_{p\le (\log |d|)^2} \Big(1 -\frac{\chi_d(p)}{p}\Big)^{-1},  
\end{equation} 
so that the extreme values of $L(1,\chi_d)$ over all $|d|\le X$ are bounded above by $(2+o(1)) e^{\gamma} \tau_{\max}$ and below by 
$(\tfrac 12+ o(1)) \zeta(2)/(e^{\gamma} \tau_{\max})$.   There is still a gap between these GRH bounds and the probabilistic conjecture in \eqref{1.2}, but now one is off only by a factor of $2$, corresponding to the expectation based on the random model that in \eqref{1.3} we only need to take the product over primes $p\le (\log |d|)$ in order to approximate $L(1,\chi_d)$.

To summarize our discussion, the values of $L(1,\chi_d)$ have an almost periodic structure in $d$, and these values may be accurately modeled by random Euler products.   The random model gives a satisfactory description of the statistical distribution of $L(1,\chi_d)$.  It also makes predictions on the largest and smallest possible values of $L(1,\chi_d)$, but there is a large gap between these predictions and our current unconditional knowledge, and even assuming GRH there is still a factor of $2$ at issue.

Similar results may be established for the distribution at the edge of the critical strip for values in other families of $L$-functions.  For example, consider the distribution of $\zeta(1+it)$, where $t$ is chosen uniformly from $[T,2T]$ with $T\to \infty$.  These values may be modeled by the random Euler product 
\begin{equation} 
\label{1.4} 
\zeta(s,{\mathbb X}) = \prod_p \Big(1- \frac{{\mathbb X}(p)}{p^s}\Big)^{-1} = \sum_{n=1}^{\infty} \frac{{\mathbb X}(n)}{n^s}, 
\end{equation} 
where the random variables ${\mathbb X}(p)$ are independent for different primes $p$, and are all chosen uniformly from the unit circle $\{ |z| =1\}$, and extended to random variables ${\mathbb X}(n)$ over all natural numbers $n$ by multiplicativity.  As before, the product and series both converge almost surely when Re$(s) >\tfrac 12$.  Then the statistical distribution of $\zeta(1+it)$ is identical to that of $\zeta(1,{\mathbb X})$ (equivalently of $\zeta(1+iy,{\mathbb X})$ for any real $y$).  We can also formulate an \emph{almost peridoicity} result:  For any $\epsilon >0$ we can find a sequence of \emph{almost periods} $\tau_n$, with $\tau_n\to \infty$ and $|\tau_{n+1} -\tau_n|$ bounded, such that  for $T$ sufficiently large (in terms of any fixed almost period $\tau$) one has $|\zeta(1+it + i\tau) -\zeta(1+it)| < \epsilon$ for almost all $t\in [T,2T]$.   The sequence of almost periods are found by requiring $p^{i\tau} \approx 1$ for all primes $p$ up to some point.   For a study of the distribution of $\zeta(1+it)$, with a focus on uniformity, see Lamzouri \cite{Lam1}.

There is an extensive literature concerned with distribution at the edge of the critical strip, and we end this section with references to some further examples.  We motivated our discussion of $L(1,\chi_d)$ with the class number formula, which (for negative fundamental discriminants) shows that $2\sqrt{|d|} L(1,\chi_d)/(2\pi)$ is quantized to be an integer.  This raises questions on the \emph{granularity} of the distribution of $L(1,\chi_d)$, and shows that in very short scales there must be arithmetic deviations from the random model.  These questions are related to the problem of understanding how many imaginary quadratic fields there are with any given class number (see \cite{HJKMP, So4, Lam5}).   For positive fundamental discriminants, the class number formula relates $L(1,\chi_d)$ to the product of the class number and the regulator which cannot in general be separated from each other.  One way to get around this problem is to order the real quadratic fields 
by the size of their regulator rather than by discriminant, and this ordering has a pleasing interpretation in terms of lengths of closed geodesics on the hyperbolic surface $PSL(2,{\mathbb Z}) \backslash {\mathbb H}$.  The study of $L(1,\chi_d)$, or the class number $h(d)$, when $d$ is ordered in 
this way was initiated by Sarnak \cite{Sar}; it is closely related to specializing discriminants $d$ in suitable quadratic sequences (for example, of the form $4n^2+1$, or $n^2+ 4$), and for recent investigations see \cite{DL, Lam6, Raulf}.     For a small sample of investigations in other families of $L$-functions, see  \cite{CogMi, Duke, XLi, Molteni, Luo}.

\section{Selberg's central limit theorem}

\noindent In the previous section we discussed the distribution of values of $L$-functions at the edge of the critical strip.  In fact, similar results hold for the value distribution inside the critical strip, but keeping to the right of the critical line.  As an illustration, consider the problem of the distribution of values of $\zeta(\sigma+it)$ where $\tfrac 12 < \sigma \le 1$ is fixed, and $t$ is chosen uniformly from $[T,2T]$ with $T \to \infty$.  The random $\zeta(s,{\mathbb X})$ defined in \eqref{1.4} still converges when Re$(s) = \sigma > \tfrac 12$, and one can show that $\zeta(\sigma+ it)$ is distributed like $\zeta(\sigma, {\mathbb X})$.  To give a very brief indication of the proof, one can show that for any parameter $1\le N\le T$ 
\begin{equation} 
\label{2.1} 
\frac 1T \int_T^{2T} \Big|\zeta(\sigma+ it) - \sum_{n\le N} \frac{1}{n^{\sigma+it}}\Big|^2 dt = O\Big( \sum_{n>N} \frac{1}{n^{2\sigma}}\Big) = O( N^{1-2\sigma}), 
\end{equation} 
which parallels 
$$ 
\E \Big[ \Big| \zeta(\sigma, {\mathbb X})  - \sum_{n\le N} \frac{{\mathbb X}(n)}{n^{\sigma}}\Big|^2 \Big]  = \sum_{n> N} \frac{1}{n^{2\sigma}} = 
O(N^{1-2\sigma}). 
$$ 
Since $\sigma>\tfrac 12$, the term $N^{1-2\sigma}$ tends to $0$ provided $N$ tends to infinity with $T$, and for such $N$ it follows that for most $t\in [T,2T]$ one has $\zeta(\sigma +it) \approx \sum_{n\le N} n^{-\sigma +it}$.   If now $N$ tends slowly to infinity with $T$, then we can show that $\sum_{n\le N} n^{-\sigma+it}$ is distributed like $\sum_{n\le N} {\mathbb X}(n)/n^{\sigma}$, by matching the moments of both quantities for example.  
This is a classical result (see Chapter XI of \cite{Tit}), and a recent quantitative study has been made in \cite{LaLeRa1}.

As with the distribution of $\zeta(1+it)$, there is an almost periodic structure in the values of $\zeta(\sigma+it)$.  The partial sums $\sum_{n\le N} n^{-\sigma-it}$ clearly have an almost periodic structure --- if $n^{i\tau} \approx 1$ for all $n\le N$, then $\tau$ will be an almost period for these partial sums --- and as we noted above $\zeta(\sigma+it)$ can often be approximated by such partial sums.  

For $\frac 12 < \sigma \le 1$, the values $\zeta(\sigma, {\mathbb X})$ are distributed densely in the complex plane; indeed, 
for any given complex number $z$ and any $\epsilon>0$, with positive probability (depending on $z$ and $\epsilon$) one has 
$|\zeta(\sigma, {\mathbb X})-z| < \epsilon$.  This is not hard to show, starting with the fact that $\log \zeta(\sigma, {\mathbb X})$ is 
essentially $\sum_p {\mathbb X}(p)/p^{\sigma}$.  It follows that the set $\{ \zeta(\sigma+it): \ t\in {\mathbb R}\}$ is dense in ${\mathbb C}$.   
A related striking \emph{universality} result of Voronin \cite{Vor} states that if $f$ is any non-vanishing continuous function in 
$|z|\le r$ with $0< r<\frac 14$, then there exist arbitrarily large values $t\in {\mathbb R}$ such that $|\zeta(\frac 34+it + z) - f(z)| < \epsilon$ for all $|z|\le r$.  
In other words, the zeta function in a disc of radius $r$ around $\frac 34+it$ can be made to mimic any given analytic function that does not take the value $0$.  
The value $0$ must be excluded in view of the Riemann Hypothesis!  There are more precise versions of this result, but we do not pursue this direction 
further, pointing instead to \cite{Bagchi, Kow, LaLeRa2} for recent related work.

\smallskip 

We now turn to the distribution of values of $\zeta(\tfrac 12 + it)$, which forms the main focus of this article.  The random Euler product $\zeta(s,{\mathbb X})$ defined in \eqref{1.4} does not converge for $s=\tfrac 12$.  Indeed, there is no almost periodic structure to the values $\zeta(\tfrac 12+it)$, and on the critical line the zeta-function cannot typically be understood simply from a knowledge of $p^{it}$ for small primes $p$.  Instead we have the following fundamental result of Selberg.

\begin{theorem}[Selberg \cite{Sel, Sel2}]  If $T$ is large, and $t$ is chosen uniformly from $[T,2T]$, then $\log \zeta(\tfrac 12+it)$ is distributed like a complex Gaussian with mean $0$ and variance $\log \log T$.  In particular, Re$(\log \zeta(\frac 12+it))$ and Im$(\log \zeta(\frac 12+it))$ are distributed like real Gaussians with mean $0$ and variance $\frac 12 \log \log T$.  
\end{theorem} 

To clarify normalizations, we recall that a standard complex Gaussian (of mean $0$ and variance $1$) has density $\frac{1}{\pi} e^{-|z|^2}$, and that its real and imaginary part are independent real Gaussians with mean $0$ and variance $\frac 12$.  Selberg's theorem gives that for any fixed box ${\mathcal B}$ in the complex plane, as $T \to \infty$ one has 
$$ 
\frac 1T \text{meas} \Big\{ T\le t\le 2T, \ \ \frac{\log \zeta(\frac 12+it)}{\sqrt{\log \log T}} \in {\mathcal B}\Big\} \to \frac 1\pi \int_{x+iy\in \mathcal B} e^{-x^2 -y^2} dx dy. 
$$  
In Selberg's theorem we may omit the countably many zeros of $\zeta(s)$ where the logarithm is not defined.  For $t$ not equalling the ordinate of a zero of $\zeta(s)$, the argument of $\zeta(\tfrac 12+it)$ (that is, Im$(\log \zeta(\tfrac 12+it))$) is defined by continuous variation along the straight lines from $2$ (where the argument is taken to be zero) to $2+it$ and thence to $1/2+ it$. 

Here is a striking illustration of the difference between the value distributions of $\zeta(\tfrac 12+it)$ and $\zeta(\sigma +it)$ for $1\ge \sigma >\frac 12$.  Typically $|\zeta(\sigma+it)|$ is of constant size, for example taking values between $1/2$ and $2$ with positive probability.  On the other hand, Selberg's theorem implies that for any fixed $V$ and large $T$ 
\begin{equation} 
\label{2.2} 
\frac 1T \text{meas} \Big\{ T\le t\le 2T, \ \ \frac{\log |\zeta(\frac 12+it)|}{\sqrt{\frac 12 \log \log T}} \ge V \Big\} \sim \frac{1}{\sqrt{2\pi}} \int_V^{\infty} 
e^{-x^2/2} dx, 
\end{equation} 
so that  $|\zeta(\tfrac 12+it)|$ is large (say $> \exp(\epsilon \sqrt{\log \log T})$) nearly half the time, or $|\zeta(\tfrac 12+it)|$ is small (below $\exp(-\epsilon \sqrt{\log \log T})$) nearly half the time.  We noted earlier that the set $\{ \zeta(\sigma+it): t\in {\mathbb R}\}$ is dense in the complex plane.   It is rare to find values of $\zeta(\frac 12+it)$ of constant size, and whether the set $\{\zeta(\tfrac 12+it): t\in {\mathbb R}\}$ is dense in ${\mathbb C}$ remains an intriguing open problem.   This question was raised first by Ramachandra; for partial progress see \cite{KoNi}.

The argument principle, together with the functional equation for $\zeta(s)$ and Stirling's formula, may be used to show that  $N(t)$, the number of zeros of $\zeta(s)$ with real part between $0$ and $1$ and imaginary part between $0$ and $t$, satisfies 
\begin{equation} 
\label{2.25}
N(t) = \frac{t}{2\pi} \log \frac{t}{2\pi} -\frac{t}{2\pi} + \frac 78 + S(t) + O\Big( \frac 1t\Big), \text{  where  }  S(t) = \frac 1\pi \text{arg} \zeta(\tfrac 12 +it).  
\end{equation}
Thus Selberg's theorem for Im$(\log \zeta(\frac 12+it))$ shows that the remainder term in the asymptotic formula for $N(t)$ has Gaussian fluctuations.

We now give a brief, oversimplified, description of the ideas behind Selberg's theorem; we caution the reader that some statements below should be taken as merely indicative, and not interpreted as being literally correct.  Taking logarithms in the Euler product for $\zeta(s)$, we may write 
$$ 
\log \zeta(s) = \sum_{p, k} \frac{1}{kp^{ks}} = \sum_{n=2}^{\infty} \frac{\Lambda(n)}{\log n} \frac{1}{n^s}, 
$$ 
where the sums above are over prime powers $p^k$, and $\Lambda(n)$ is the von Mangoldt function which equals $\log p$ if $n=p^k$ and $0$ otherwise.  The series above converges absolutely when Re$(s)>1$, and it certainly does not converge on the critical line Re$(s) =\frac 12$.   Nevertheless, we might hope that a truncated sum over prime powers might serve as an approximation to $\log \zeta(s)$ (thinking of $s=\frac 12+it$ with $T\le t\le 2T$).  This forms the first step in Selberg's argument, who finds an expression of the form 
\begin{equation} 
\label{2.3} 
\log \zeta(s) = \sum_{2 \le n\le x} \frac{\Lambda(n)}{n^s \log n} + Z_x(s), 
\end{equation} 
where $Z_x(s)$ is a remainder term that may be thought of as the contribution from zeros $\rho$ of $\zeta(s)$ with $|\rho -s| \le 1/\log x$.  By a complicated argument Selberg showed how the sum over zeros may in turn also be bounded in terms of sums over primes, and thus shown to be small on average.  An alternative argument of Bombieri and Hejhal  \cite{BomHej}  avoids some of Selberg's difficulties by bounding the average values of $Z_x(s)$ instead of seeking point-wise bounds.  Nevertheless, these arguments are technically involved; they are simpler if the Riemann hypothesis is assumed, but can be established unconditionally by relying on a subtle zero-density estimate for zeros of $\zeta(s)$ near the critical line (established by Selberg).  Although we haven't made the relation \eqref{2.3} precise, we give a couple of remarks that may be helpful in thinking about such relations.  Firstly, one can think of such relations as variants of the explicit formula connecting zeros and primes.  Secondly, in addition to the Euler product, the zeta function possesses a Hadamard product over its zeros 
\begin{equation} 
\label{2.4} 
s(s-1) \pi^{-s/2} \Gamma(s/2) \zeta(s) = e^{Bs} \prod_{\rho} \Big(1 -\frac s{\rho}\Big) e^{s/\rho}, 
\end{equation} 
where the product is over all non-trivial zeros of the zeta-function, and $B$ is a constant.  The relation \eqref{2.3} has the flavor of a hybrid Euler--Hadamard product (see \cite{GHK} for work in this direction), taking some primes and some zeros, and it is natural to expect an inverse relationship (or uncertainty principle) between the number of primes that one must take versus the number of zeros that are needed.   

Returning to the argument, in the range $x\le T$, the remainder term $Z_x(s)$ in \eqref{2.3} is typically of size $O(\log T/\log x)$ --- this corresponds to the expected number of zeros of $\zeta(s)$ within $1/\log x$ of $\frac 12 +it$.  If we choose $x= T^{1/(\log \log T)^{\frac 14}}$ for example, then $\log T/\log x  = (\log \log T)^{\frac 14}$ is small in comparison to the typical expected size of $\log \zeta(s)$, which is $\sqrt{\log \log T}$, and therefore the remainder term is negligible.  In other words, with this choice of $x$, the proof of Selberg's theorem reduces to establishing the Gaussian nature of 
\begin{equation} 
\label{2.5} 
\sum_{2\le n\le x} \frac{\Lambda(n)}{\log n} \frac{1}{n^s}  = \sum_{p\le x} \frac{1}{p^s}  + \frac 12 \sum_{p\le \sqrt{x}} \frac{1}{p^{2s}} +
 \sum_{\substack{ p^k \le x\\ k\ge 3}} \frac{1}{kp^{ks}}. 
 \end{equation} 
 The contribution from prime powers $p^k$ with $k\ge 3$ is $O(1)$ and may be omitted.  The contribution from the squares of primes is also negligible; it is $\frac 12 \sum_{p\le \sqrt{x}} 1/p^{1+2it}$ which behaves roughly like $\frac 12 \log \zeta(1+2it)$ and so is of constant size typically.
We are left with the contribution of just the primes, which we may understand by computing moments.  If $k$ and $\ell$ are any natural numbers 
then, for large $T$,  
 \begin{equation} 
 \label{2.6} 
 \frac 1T \int_T^{2T} \Big(\sum_{p\le x} \frac{1}{p^{1/2+it}}\Big)^k \Big( \sum_{p\le x} \frac{1}{p^{1/2-it}}\Big)^{\ell} dt =
 \begin{cases} 
  (1+o(1)) k! (\log \log T)^{k} &\text{ if } k =\ell\\ 
  o(T) &\text{ if  } k \neq \ell.  
 \end{cases} 
 \end{equation}
 These moments match asymptotically the moments of a complex Gaussian with mean $0$ and variance $\log \log T$, from which Selberg's theorem would follow.

 To give a justification for \eqref{2.6}, we discuss an orthogonality relation for Dirichlet polynomials, which we shall find useful in the sequel.  Roughly speaking, integrals over $[T,2T]$ may be thought of as possessing $T$ ``harmonics'' that can distinguish between the functions $f_n(t) = n^{it}$ for  
 natural numbers $n$ going up to about $T$.  More precisely, suppose $\Phi$ is a smooth function approximating the indicator function of $[1,2]$.  Then, if $\max(M,N) \le T/\log T$,  
 \begin{align} 
 \label{2.7}
 \int \sum_{m \le M} a(m) m^{it}  \overline{\sum_{n\le N} b(n) n^{it}} \Phi \Big( \frac{t}{T} \Big) dt 
 &= \sum_{m=n} a(m)\overline{b(n)} T {\hat \Phi}(0)  + \sum_{m\neq  n} a(m) \overline{b(n)} T {\hat \Phi}\Big( T \log \frac{n}m\Big) \nonumber \\ 
&\sim T {\hat \Phi}(0)\sum_{m=n} a(m)\overline{b(n)},
\end{align}
where the contribution of the ``off-diagonal'' terms $m\neq n$ is negligible because $T |\log (m/n)| \gg T |m-n|/|m+n| \ge T/(M+N)$ is large and the 
Fourier transform $\hat \Phi$ decays rapidly.   

Write $(\sum_{p\le x} 1/p^{1/2+it})^k =\sum_{n\le x^k} a_k(n)/n^{1/2+it}$, so that $a_k(n) =0$ unless $n$ has exactly $k$ prime factors.  If $n$ has prime factorization $p_1^{e_1} \cdots p_r^{e_r}$ with $e_1+\ldots +e_r = k$ then $a_k(n) = k!/(e_1! \cdots e_r!)$.  Then an application of \eqref{2.7} 
shows that the moment in \eqref{2.6} is 
$$ 
\sim  T \sum_{m=n \le x^k} \frac{a_k(n) a_\ell(n)}{n} .
$$ 
If $k\neq \ell$ then either $a_k(n)$ or $a_\ell(n)$ must be zero, and this case of \eqref{2.6} follows.  If $k=\ell$, then the diagonal terms are dominated by integers with $k$ distinct prime factors, and so the above is 
$$ 
\sim T k! \sum_{n\le x^k} \frac{a_k(n)}{n} = T k! \Big( \sum_{p\le x} \frac 1p \Big)^k \sim T k! (\log \log x)^k,
$$ 
and since $\log \log x$ and $\log \log T$ are close, the other case in \eqref{2.6} follows.

This concludes our sketch of the ideas behind Selberg's theorem.  Two alternative approaches that work for $\log |\zeta(\tfrac 12+it)|$ are given in 
\cite{Laur, RaSo3}.   These avoid the subtle zero density estimates near the critical line, and it would be of interest to extend such approaches to Im$(\log \zeta(\frac 12+it))$.

%%Applications --- relevance to moments; zeros of linear combinations ? ?

\section{Analogues of Selberg's theorem in families of L-functions} 

\noindent Selberg's theorem discussed above applies not only to the Riemann zeta-function, but more generally to a large class of $L$-functions.  For example, in \cite{Sel} Selberg introduced what is now known as the \emph{Selberg class} of $L$-functions, which formalizes some of the observed properties of automorphic $L$-functions and is expected to coincide with this class.  For a primitive  $L$-function in the Selberg class (or, if one prefers, for a cuspidal automorphic $L$-function for $GL_n({\mathbb Q})$), one expects that $\log L(\frac 12+it)$ with $T \le t\le 2T$ is distributed like a complex Gaussian with mean $0$ and variance $\log \log T$.  The key ingredient needed to make this precise is an analogue of the zero density estimate close to the critical line, and this is known for $GL_1$ and $GL_2$; in the general case, GRH must be assumed (see \cite{Sel, BomHej} for more details).

Interesting differences arise when we consider analogues of Selberg's theorem for central values in families of $L$-functions.  There are three categories into which families of $L$-functions fall, and we illustrate these with examples.  Unlike Selberg's \emph{Theorem}, the analogous central limit theorems that we formulate in these families are still \emph{conjectural}, and these conjectures were first formulated by Keating and Snaith \cite{KeSn1}.

\emph{Unitary families.}  A typical example is the family of all Dirichlet characters $\chi \pmod q$, with $q$ a large prime (for simplicity).   The question is to understand the distribution of $\log L(\tfrac 12,  \chi)$ as $\chi$ ranges over all primitive characters $\chi \pmod q$ (if $q$ is prime, this is equivalent to $\chi$ not being the principal character).   We must discard potential characters with $L(\frac 12,\chi)=0$, but in fact it is conjectured that $L(\tfrac 12, \chi)\neq 0$ for all Dirichlet $L$-functions.   This situation is expected to be exactly as in Selberg's theorem, and the Keating--Snaith conjecture for this family states that for large $q$ the distribution of $\log L(\tfrac  12, \chi)$ is approximately a complex Gaussian with mean $0$ and variance $\log \log q$. 
% --- say something about the argument here ---  
In particular $\log |L(\frac 12, \chi)|$ is (conjecturally) distributed like a real Gaussian  with mean $0$ and variance $\frac 12 \log \log q$, so that (like $|\zeta(\tfrac 12+it)|$) roughly half the time $|L(\tfrac 12, \chi)|$ is as large as $\exp(\epsilon \sqrt{\log \log q})$ and the other half of the time it is as small as $\exp(-\epsilon \sqrt{\log \log q})$.

Another example of this type is the family of twists by Dirichlet characters of a fixed newform $f$.  
The family $\zeta(\frac 12+it)$ with $T\le t\le 2T$ may also be thought of as an example of a unitary family.

\emph{Symplectic families.}   Consider the family of quadratic Dirichlet $L$-functions $L(s,\chi_d)$, where $d$ ranges over fundamental discriminants with $|d|\le X$.  
The values $L(\tfrac 12,\chi_d)$ are real, and GRH predicts that they are all non-negative (else there would be a real zero of $L(s,\chi_d)$ between $1/2$ and $1$).  Further, the values $L(\frac 12, \chi_d)$ are all expected to be non-zero (a conjecture of Chowla, which is a special case of the belief that $L(\frac 12 , \chi)\neq 0$ for all Dirichlet characters $\chi$).   The Keating--Snaith conjecture for this family predicts that the values $\log L(\frac 12, \chi_d)$ are distributed like a real Gaussian with mean $\frac 12 \log \log X$ and variance $\log \log X$.   Since the mean is positive, the values of $L(\tfrac 12, \chi_d)$ are (conjecturally) of typical size $(\log X)^{\frac 12+ o(1)}$.

\emph{Orthogonal families.}  These families arise naturally in the context of modular forms, and we give a couple of prototypical examples.  Let $k$ be an even integer, and consider the family ${\mathcal H}_k$ of all weight $k$ modular forms for the full modular group $SL_2({\mathbb Z})$ that are also eigenfunctions of all Hecke operators.  Associated to such a form $f$ is its $L$-function, which we normalize so that the functional equation connects values at $s$ to $1-s$:
$$ 
\Lambda(s,f) = (2\pi)^{-s} \Gamma(s+ \tfrac {k-1}{2}) L(s,f) = i^{k} \Lambda(1-s, f). 
$$
In the case $k\equiv 2 \pmod 4$, the sign of this functional equation is $-1$, and all the central values $L(\tfrac 12, f)$ are zero.  In the case $k\equiv 0\pmod 4$, the sign of the functional equation is $+1$, and we ask for the distribution of $L(\tfrac 12, f)$ (or, in keeping with Selberg's theorem, $\log L(\tfrac 12,f)$).  In this situation, a remarkable result of Waldspurger \cite{W} (see also \cite{KZ} for an explicit version) relates these central $L$-values to the squares of Fourier coefficients of a half-integer weight modular form associated to $f$ (namely its Shimura correspondent).  As a byproduct, we know that $L(\frac 12, f)$ is non-negative, and it is conjectured never to be zero.  The Keating--Snaith conjectures predict that for large $k \equiv 0 \pmod 4$, the values $\log L(\tfrac 12, f)$ are distributed like a real Gaussian with mean $-\frac 12 \log \log k$ and variance $\log \log k$.   Since the mean is negative, the values $L(\frac 12, f)$ in this family are typically small, of size $(\log k)^{-\frac 12+ o(1)}$.   

A related example is to fix a newform $f$, and to consider the family of quadratic twists of $f$.  Once again normalizing so that the functional equation connects $s$ and $1-s$, our interest is in the central values $L(\tfrac 12, f \times \chi_d)$, where $d$ runs over fundamental discriminants $|d|\le X$ with $d$ coprime to the level of $f$ for simplicity.   As in the previous example, half of these twists will have a functional equation with $-$ sign (where the central $L$-value vanishes), and we restrict attention to the complementary case when the sign is $+$.  Again Waldspurger's formula shows that the central $L$-values are non-negative, but it is possible for these values to be $0$.  For example, if $f$ corresponds to an elliptic curve, then the Birch--Swinnerton-Dyer conjectures predict that the central value is zero when the quadratic twist of this elliptic curve has positive rank (and the rank must also be even when the sign of the functional equation is $+$).  However, one expects that typically $L(\frac 12, f\times \chi_d) \neq 0$, and the Keating--Snaith conjectures predict further that the distribution of $\log L(\tfrac 12, f\times \chi_d)$ (where $|d|\le X$ is coprime to the level of $f$ and the twist has $+$ sign of the functional equation)  is that of a real Gaussian with mean $-\frac 12 \log \log X$ and variance $\log \log X$.

 \smallskip 
 
 The classification of families into unitary, symplectic, and orthogonal is based on the philosophy of Katz and Sarnak \cite{KaSa} which connects (conjecturally) the distribution of low lying zeros in these families to the distribution of eigenvalues near $1$ of large random matrices chosen from the corresponding classical groups --- we shall discuss these links to random matrix theory later.   We now give heuristic reasons to explain the three different Keating--Snaith conjectures, point out the obstructions to making these precise, and describe the partial progress that has been made.

 Recall that in \eqref{2.3}  we considered approximations to $\log \zeta(\tfrac 12+it)$ by Dirichlet series over prime powers of a flexible length $x$.   In \eqref{2.5} we saw that for $\zeta(\tfrac 12+it)$, the contribution of prime powers $p^k$ with $k\ge 3$ is bounded, and the contribution from prime squares is also typically small.   Finally the distribution of the sums over primes could be understood by computing moments.   We now consider analogues of this calculation for the families discussed above, and the key difference in the orthogonal and symplectic cases will arise in the contribution of squares of primes.    
 
 Let us first look at the unitary family of Dirichlet characters $(\bmod \, q)$ with $q$ a large prime.  Suppose that we have an approximation of the form 
 \begin{equation} 
 \label{3.1} 
 \log L(\tfrac 12, \chi) \approx \sum_{n\le x} \frac{\Lambda(n)}{\sqrt{n} \log n} \chi(n) = \sum_{p\le x} \frac{\chi(p)}{\sqrt{p}} + \frac 12 \sum_{p\le \sqrt{x}} \frac{\chi(p)^2}{p} + O(1).
 \end{equation} 
A typical character $\chi \pmod q$ is not quadratic; $\chi^2$ is then a non-principal character and the sum over prime squares above is typically of bounded size, behaving a lot like $\log L(1, \chi^2)$.  We are left with the sum over primes, and if $x$ is a small power of $q$, then we can understand the moments of this sum 
(much as in \eqref{2.6}) using the orthogonality relation for the characters $\pmod q$ (in place of \eqref{2.7}).   This gives a heuristic justification for the Keating--Snaith conjectures in this family, and the missing ingredient is the very first step which may fail badly, for example,  if $L(\frac 12, \chi) =0$ for many characters $\chi \pmod q$.  
 
 Consider next the symplectic example of quadratic Dirichlet $L$-functions $L(s,\chi_d)$ with $d$ ranging over fundamental discriminants $|d|\le X$.  Suppose that an approximation as in \eqref{3.1} holds.  Since $\chi_d$ is a quadratic character, note that the squares of primes in \eqref{3.1} have $\chi_d(p)^2=1$ (ignoring the primes $p$ that divide $d$), and so these terms contribute 
 $$ 
 \frac 12 \sum_{p\le \sqrt{x}} \frac{1}{p} \sim \frac 12 \log \log x \sim \frac 12 \log \log X, 
 $$ 
 if $x$ is a small power of $X$.   Thus the prime square terms account for the mean of  $\log L(\frac 12, \chi_d)$ being $\sim \frac 12 \log \log X$ in the Keating--Snaith conjectures.  If $x$ is a small power of $X$, we may compute the moments of the sum over primes:
 $$ 
 \sum_{|d|\le X} \Big( \sum_{p\le x} \frac{\chi_d(p)}{\sqrt{p}}\Big)^k  = \sum_{p_1, \ldots, p_k \le x} \frac{1}{\sqrt{p_1\cdots p_k}} 
 \sum_{|d| \le X} \Big( \frac{d}{p_1\cdots p_k}\Big). 
 $$  
 The inner sum over $d$ may be viewed as a character sum $(\bmod\, {p_1\cdots p_k})$.  This character is principal if $p_1\cdots p_k$ is a square, and we get a main term here, while if $p_1\cdots p_k$ is not a square we may expect the character sum to cancel out (and this can be justified if $x^k$ is small in comparison to $X$).  
 The product $p_1\cdots p_k$ can be a square only if $k$ is even, and the primes $p_1$, $\ldots$, $p_k$ can be paired off into $k/2$ equal pairs.  With a little calculation, this shows that the moments of the sum over primes match the moments of a real Gaussian with mean $0$ and variance $\sum_{p\le x} 1/p \sim \log \log X$.  Taking into account the shift in mean arising from the prime square terms, this gives a heuristic justification for the Keating--Snaith conjecture.

 Finally let us look at  the orthogonal family of quadratic twists of a newform in the case where the sign of the functional equation is $+$.  The $L$-function $L(s,f\times  \chi_d)$ is given by an Euler product, the $p$-th factor of which (for a prime $p$ not dividing  the level of the form) takes the  shape 
 $$
 \Big(1-\frac{ \alpha_p \chi_d(p)}{p^s}\Big)^{-1} \Big (1-\frac{\beta_p\chi_d(p)}{p^s}\Big)^{-1}, 
 $$ 
where $\alpha_p \beta_p=1$ and $\alpha_p + \beta_p = \lambda(p)$ is the normalized Hecke eigenvalue of $f$ (normalized so that the Deligne bound gives $|\lambda(p)| \le 2$).   The logarithm of this Euler factor is  
$$
\sum_{k=1}^{\infty} (\alpha_p^k + \beta_p^k) \frac{\chi_d(p^k)}{kp^{ks}},
$$ 
 and  in  analogy with \eqref{2.3}, \eqref{2.5}, \eqref{3.1},  we may hope to approximate $\log L(\frac 12, f\times \chi_d)$ by 
\begin{align}
\label{3.2} 
&\sum_{p\le x} \frac{(\alpha_p + \beta_p)\chi_d(p)}{\sqrt{p}} + \frac 12 \sum_{p\le \sqrt{x}}\frac{(\alpha_p^2 + \beta_p^2)\chi_d(p)^2}{p}  + O(1) 
\nonumber \\
= &\sum_{p\le x} \frac{\lambda(p)\chi_d(p)}{\sqrt{p}}  + \frac 12 \sum_{p\le \sqrt{x}} \frac{\lambda(p)^2 -2}{p} + O(1). 
\end{align}
If the discriminants $d$ go up to size $X$, and $x$ is a small power of $X$, then the distribution of $\sum_{p\le x } \lambda(p) \chi_d(p)/\sqrt{p}$ may be determined  by  computing moments (similarly to the discussion for $L(\frac 12,\chi_d)$).   The prime terms in \eqref{3.2} are distributed like a real Gaussian with mean $0$ and variance 
\begin{equation} 
\label{3.3} 
\sum_{p\le x} \frac{\lambda(p)^2}{p} \sim \log \log x \sim \log \log X,
\end{equation}
by Rankin--Selberg theory.  In view of \eqref{3.3}, the prime square terms in \eqref{3.2} contribute 
$$ 
\frac 12 \sum_{p\le \sqrt{x}}  \frac{\lambda(p)^2 -2}{p} \sim - \frac 12 \log \log \sqrt{x}  \sim -\frac 12 \log \log X.
$$ 
This justifies the Keating--Snaith conjecture for this family.  

\smallskip 

In all these heuristics, it is the first step of connecting $\log L(\tfrac 12)$ to sums over prime powers that is a serious stumbling-block.   Indeed if $L(\frac 12)$ is 
zero (or if there is a zero very close to $\frac 12$) for many elements in the family, then the Keating--Snaith conjectures would not hold.   This problem does  not arise in the continuous Selberg theorem, since the points $t$ with $\frac 12+it$ very close to a zero of $\zeta(s)$ have small measure and thus do not affect the distribution.

The problem of non-vanishing of $L$-functions has been investigated extensively, but in general it remains a challenge to show that almost all elements in a family have non-zero central value.  More often, progress towards this problem focusses on showing that a positive proportion of $L$-functions in a family have non-zero central value.   To give a few examples: in the family of  Dirichlet characters $\chi \pmod q$, Khan and Ngo \cite{KhanNgo} have shown that at least $\frac 38$ of these characters have $L(\frac 12, \chi) \neq 0$;   in the family of quadratic Dirichlet $L$-functions, Soundararajan \cite{So5} shows that a proportion at least $\frac 78$ of 
such central values are non-zero; in the family ${\mathcal H}_k$ of all Hecke eigenforms of weight $k \equiv 0\pmod 4$ for the full modular group, with $k\le K$, Iwaniec and Sarnak \cite{IwSa} show that at least $\frac 12$ of the central values are non-zero, and improving this proportion (in a certain sense) would have consequences for the existence of Landau--Siegel zeros of Dirichlet $L$-functions.   

There are some situations where, for deep algebraic reasons, one can show that most central values in a family are non-zero, but these arguments do not appear to control the size of the central value, or to deal with the possibility that there might be a zero very near $\frac 12$.  For example, Chinta 
\cite{Chinta} (following work of Rohrlich \cite{Rohr}) has shown that if $E$ is an elliptic curve over ${\mathbb Q}$ then for all but $O(q^{\frac 78})$ of the 
Dirichlet characters $\bmod \, q$ (with $q$ a large prime) one has $L(\tfrac 12, E\times \chi) \neq 0$.  This exploits the fact (established by Shimura) that if $\chi^{\sigma}$ is a Galois conjugate of the character $\chi$, then the vanishing of $L(\frac 12,  E\times \chi)$ is equivalent to the vanishing of $L(\frac 12, E\times \chi^{\sigma})$ 
(the algebraic parts of these $L$-values are Galois conjugate).    Another example where algebraic techniques are very successful concerns the family of quadratic twists  of an elliptic curve.  In special cases, Smith \cite{Smith} has shown that the (algebraic) rank of quadratic twists of elliptic curves  is typically $0$ (when the sign of the functional equation is $+$) or $1$ (when the sign is $-$).  The Birch--Swinnerton-Dyer conjecture (on which there 
has been a lot of progress in the cases of rank $0$ and $1$) would then yield Goldfeld's conjecture that the central $L$-values are typically non-zero when the sign of the functional equation is $+$.

If there is a zero at or very near $\frac 12$, we might expect that its effect is to make $|L(\frac 12)|$ unusually small.   This observation was made in Soundararajan \cite{So}, where it was shown (assuming GRH) that $\log |L(\frac 12)|$ can be bounded from above using Dirichlet series over prime powers of flexible length; we shall discuss this in more detail in \S 6.   It was also observed in \cite{So} that one could (assuming a suitable GRH) establish a one sided version of the Keating--Snaith conjecture, showing that the frequency with which $\log |L(\frac 12)| \ge \text{Mean}  +\lambda \sqrt{\text{Var}}$ is bounded above 
by the expected Gaussian $\frac{1}{\sqrt{2\pi}} \int_{\lambda}^{\infty} e^{-x^2/2} dx$; here $\lambda$ is a fixed real number, the size of the family is assumed to grow.  
Further, if one knew that most elements in the family did not have a zero near $\frac 12$ (which, for example, would follow from the ``one level density'' conjectures in Katz and Sarnak \cite{KaSa}) then the Keating--Snaith conjecture for $\log |L(\frac 12)|$ would follow. 

 Such one sided central limit theorems were first made precise (and unconditional) by Hough \cite{Hough} in certain families of $L$-functions.   Hough's approach relies on knowledge of a zero density estimate putting most low lying zeros of $L$-functions in the family close to the critical line --- an analogue of Selberg's zero density estimate for the zeta function, mentioned in \S 2.   For example, Hough's approach would work for $\log |L(\frac 12, \chi)|$ in the unitary family of Dirichlet characters $\chi \pmod  q$, or $\log |L(\frac 12, \chi_d)|$ in the symplectic family of quadratic Dirichlet $L$-functions, or in the orthogonal family  $\log L(\frac 12, f)$ 
 where $f$ ranges over Hecke eigenforms of weight $k \equiv 0 \pmod 4$ for the full modular group.

An alternative approach to this half of the Keating--Snaith conjectures is developed in Radziwi{\l \l} and Soundararajan \cite{RaSo2}.  This method is arguably simpler and also more widely applicable, relying only on knowledge of the first moment ``$+$ epsilon'' in the family, and avoiding zero density estimates (which require knowledge of the second moment ``$+$ epsilon'').   In \cite{RaSo2} the method is illustrated for the family of quadratic twists of an elliptic curve (with positive sign of 
the functional equation), where the zero density estimates required in Hough's approach are not known.   Conjecturally the central values in this family  (when non-zero) 
measure (after accounting for quantities such as Tamagawa factors that are relatively easy to understand) the size of the Tate--Shafarevich group for the twisted elliptic curve.  The Keating--Snaith conjecture thus predicts that the sizes of Tate--Shafarevich groups in the family of quadratic twists have a log normal distribution, 
with prescribed means and variance  (see Conjecture 1 in \cite{RaSo2}).   The method applies to quadratic twists of any newform (holomorphic or Maass form), and thus  (by Waldspurger's formula) also gives information on the size of Fourier coefficients of half-integer weight modular forms, establishing that these are typically a little bit smaller than the conjectured Ramanujan bounds.

Another application where this method works is to the problem of the fluctuations of a quantum observable for the modular surface.   Let $\psi$ denote a fixed even Hecke-Maass form for the full modular group, and let $\phi_j$ denote an even Hecke-Maass form with eigenvalue $\lambda_j$.  The problem is to understand $\mu_j(\psi)= \int_{PSL_2({\mathbb Z})\backslash {\mathbb H}} \psi(z) |\phi_j(z)|^2 \frac{dx dy}{y^2}$ for large eigenvalue $\lambda_j$.    For generic  hyperbolic surfaces, it has been suggested in the physics literature \cite{Eckhardtetal}  that similar quantum fluctuations have a Gaussian distribution.  In the case of the modular group, $|\mu_j(\psi)|^2$ is related to the central value $L(\frac 12, \psi \times \phi_j \times \phi_j)$, so that the Keating--Snaith conjectures predict that it is in fact $\log |\mu_j(\psi)|$ (rather than $\mu_j(\psi)$ itself) that has a normal distribution.   A one sided central limit theorem for $\log |\mu_j(\psi)|$ is obtained in Siu \cite{Siu}, and in particular it follows that $\lambda_j^{\frac 14} |\mu_j(\psi)| = o(1)$ for almost all eigenfunctions $\phi_j$.  

We have already discussed that the problem of non-vanishing of central $L$-values is a barrier to obtaining   lower bounds towards the Keating--Snaith conjectures.  There are two analytic techniques that produce  a positive proportion of non-zero central values of $L$-functions in families: (i) the mollifier method, which is unconditional and relies on knowledge of two moments (``$+$ epsilon'') and (ii) understanding $1$-level densities of low lying zeros, which is conditional on GRH and is not always guaranteed to yield a non-zero proportion.  Both of these methods may be refined to permit an understanding of the typical size of non-zero $L$-values that are produced \cite{SoundMFO}.  Here are two such sample results.  In the family of quadratic Dirichlet $L$-functions, where we know \cite{So5} that $\frac 78$ of the fundamental discriminants $|d|\le X$ satisfy $L(\tfrac 12,\chi_d) \neq 0$,  we may establish that for any interval $(\alpha,\beta)$ of ${\mathbb R}$ and large $X$ 
$$ 
\#\Big\{ |d|\le X: \ \ \frac{\log |L(\tfrac 12,\chi_d)| - \tfrac 12\log \log X}{\sqrt{\log \log X}} \in (\alpha,\beta)\Big\} 
\ge \Big( \frac 78 \frac 1{\sqrt{2\pi}}\int_{\alpha}^{\beta} e^{-x^2/2}dx +o(1)\Big) \# \{ |d|\le X\}. 
$$ 
In the family of quadratic twists of a fixed newform $f$ with positive sign of the functional equation, on GRH it is known that a proportion $\ge \frac 14$ of 
such $L$-values are non-zero (see \cite{HB7}), and we may refine this to yield (with ${\mathcal E}(X)$ denoting the set of fundamental discriminants $|d|\le X$ 
with the quadratic twist of $f$ has positive sign) 
$$ 
\#\Big\{ d\in {\mathcal E}(X): \ \ \frac{\log L(\tfrac 12, f\times \chi_d) + \frac 12 \log \log X}{\sqrt{\log \log X}}\in (\alpha,\beta)\Big\} 
\ge \Big(\frac 14 \frac{1}{\sqrt{2\pi}} \int_{\alpha}^{\beta} e^{-x^2/2}dx  + o(1)\Big) |{\mathcal E}(X)|. 
$$
Finally, we mention recent work of Bui et al  \cite{Buietal} which considers a variant of the Keating--Snaith conjectures when $L$-values 
are counted with suitable weights (which depend on ``mollified $L$-values'').  

\section{Moments of zeta and $L$-functions} 

\noindent A classical problem, going back to Hardy and Littlewood, asks for an understanding of the moments of $\zeta(\frac 12+it)$:
\begin{equation} 
\label{4.1}
M_k(T) = \int_0^T |\zeta(\tfrac 12+it)|^{2k} dt,  
\end{equation}
where $k$ is a natural number.   Hardy and Littlewood established that 
$M_1(T) \sim T \log T$ (see \cite{Tit}), and this was later refined by Ingham who showed that 
\begin{equation} 
\label{4.2} 
M_1(T) = \int_0^{T} |\zeta(\tfrac 12+it)|^2 dt = T \log \frac{T}{2\pi}  +(2\gamma -1) T + E(T), 
\end{equation}
with $E(T)=O(T^{\frac 12} \log T)$, with a further refinement in Balasubramanian \cite{Balu1} yielding $E(T) = O(T^{\frac 13+\epsilon})$.   Ingham also established an 
asymptotic for the fourth moment: $M_2(T) \sim \frac{1}{2\pi^2} T(\log T)^4$, which was refined by Heath-Brown \cite{HB1} to 
\begin{equation} 
\label{4.3} 
M_2(T) = \int_0^{T} |\zeta(\tfrac 12+it)|^4 dt = T P_4(\log T)  + O(T^{\frac 78 +\epsilon}), 
\end{equation} 
for a polynomial $P_4$ of degree $4$ with leading coefficient $1/(2\pi^2)$.
  
Despite much effort, these remain the only two cases in which an asymptotic 
formula for $M_k(T)$ is known.   To explain why, we recall that Hardy and Littlewood gave an ``approximate functional equation'' (in fact Riemann's unpublished notes had a more precise version, known now as the Riemann--Siegel formula) 
\begin{equation} 
\label{4.4} 
\zeta\Big(\frac 12+it \Big) \approx \sum_{n\le \sqrt{|t|/2\pi}} \frac{1}{n^{\frac 12+it}} + e^{i\vartheta(t)} \sum_{n\le \sqrt{|t|/2\pi}} \frac{1}{n^{\frac 12-it}},
\end{equation} 
where $e^{i \vartheta(t)} = \pi^{it/2} \Gamma((\frac 12-it)/2)/(\pi^{-it/2} \Gamma((\frac 12 +it)/2))$ is the ratio of $\Gamma$-factors in the functional equation for $\zeta(s)$.   Thus $\zeta(\frac 12+it)$ can be approximated by two Dirichlet polynomials of length about $\sqrt{|t|}$.   We saw in \eqref{2.7} that the mean square of Dirichlet polynomials of length up to $T$ could be evaluated, with the diagonal terms making the dominant contribution.   This permits the evaluation of the second moment \eqref{4.2} with Ingham's bound on the remainder term $E(T)$ (we have not discussed the cross terms that arise in squaring \eqref{4.4} but these turn out to be negligible).   Similarly, we can approximate $\zeta(\frac 12+it)^2$ by two Dirichlet polynomials of length about $|t|/2\pi$, and this leads to Ingham's asymptotic for $M_2(T)$, although the more precise form in \eqref{4.3} requires further ideas.   When $k\ge 3$, the complexity of $\zeta(\frac 12+it)^{k}$ becomes too great; to approximate it we require Dirichlet polynomials of length about $|t|^{k/2}$ (which is now larger than $|t|$), and \eqref{2.7} is no longer sufficient to estimate the mean-square of such long Dirichlet polynomials.    Let $d_k(n)$ denote the $k$-divisor function, which arises as the Dirichlet series coefficients of $\zeta(s)^k = \sum_{n=1}^{\infty} d_k(n)/n^s$ (valid for Re$(s)>1$).   One new problem that arises when considering higher moments involves the correlations 
\begin{equation} 
\label{4.5} 
\sum_{n\le x} d_k(n) d_k(n+h). 
\end{equation} 
One would like asymptotics for such quantities, uniformly in a range for $h$, and while this problem has been solved for $k=2$ (and underlies the precise asymptotics given in \eqref{4.3}), when $k=3$ or larger, asymptotics for the quantity in \eqref{4.5} remain unknown (even in the case $h=1$).

Indeed until the late 90's it was not even clear what the conjectural asymptotics for $M_k(T)$ should be.  
However in the last twenty five years, much progress has been made in understanding conjecturally 
the nature of these moments, obtaining lower bounds of the correct conjectured value (for all positive real $k$), and obtaining complementary upper bounds 
of the correct order conditional on the Riemann Hypothesis.  Similar progress has been made for moments in a number of different families of $L$-functions.   We shall discuss these conjectures and the progress towards them in the following sections, but first give some motivation for considering such moments.

One motivation for considering the moments of $\zeta(s)$ is that they capture information about the large values of $|\zeta(\tfrac 12+it)|$.   
The Lindel{\" o}f hypothesis that $|\zeta(\tfrac 12+it)| \ll_{\epsilon} (1+ |t|)^{\epsilon}$ (which is a consequence of RH) is equivalent to 
the bound $M_k(T) \ll_{k,\epsilon} T^{1+\epsilon}$ for all $k\in {\mathbb N}$.   From the approximate functional equation \eqref{4.4} it follows that 
$|\zeta(\frac 12+it)| \ll (1+|t|)^{\frac 14}$, a bound known as the \emph{convexity bound}.   Going beyond the convexity bound involves showing cancellation in the 
exponential sums in \eqref{4.4}, and has remained an active problem from its initiation by Weyl, and Hardy and Littlewood who showed early on that $|\zeta(\frac 12+ it)| \ll (1+|t|)^{\frac 16+\epsilon}$ (see \cite{Tit} and the best current exponent may be found in  \cite{Bour}).    Sharp moment estimates encode Lindel{\" o}f bounds on average, and in some cases can also yield pointwise subconvexity estimates.   For example, we note that 
Ingham's bound $E(T) \ll T^{\frac 12} \log T$ (for the error term in the second moment \eqref{4.2}) implies that 
$\int_T^{T+1} |\zeta(\tfrac 12+it)|^2 dt \ll T^{\frac 12}\log T$ from which the convexity bound $|\zeta(\tfrac 12+it)| 
\ll |t|^{\frac 14+ \epsilon}$ may be deduced.   Similarly Balasubramanian's improved estimate for $E(T)$ implies the Hardy-Littlewood-Weyl subconvexity bound $|\zeta(\frac 12+it)| \ll (1+|t|)^{\frac 16+\epsilon}$.     Similarly, Ingham's asymptotic for the fourth moment yields the convexity bound, while the more precise result \eqref{4.3} of Heath-Brown 
gives a subconvexity bound for $\zeta(s)$.  
  As a third 
example of bounds for moments that encode good point-wise bounds, we mention Heath-Brown's 
\cite{HB2} estimate for the twelfth moment 
$$ 
\int_0^{T} |\zeta(\tfrac 12+it)|^{12} dt \ll T^{2+\epsilon}, 
$$ 
which again contains the bound $|\zeta(\tfrac 12+it)| \ll |t|^{\frac 16+\epsilon}$.     

Ingham's work on the fourth moment of $\zeta(\frac 12+it)$ is also crucial in establishing ``zero density estimates'' which are bounds for the number of 
potential exceptions to the Riemann hypothesis.  These have arithmetic applications, for example playing a key role in showing that the prime number theorem 
holds in short intervals: $\pi (x+h) - \pi(x) \sim h/\log x$ provided $x^{\frac 7{12}+\epsilon} < h \le x$.  A sharp bound for the sixth moment (for instance) would lead to improvements in zero density results and in the application to the prime number theorem.   We refer to Chapter 10 of \cite{IK} for a discussion of these themes.

There is a large body of work studying analogous problems for moments of central values in families of $L$-functions, and in many cases asymptotics 
for small moments are known.   We give a few examples here, and discuss some more in \S 6.  Two motivations for studying such questions are (i) the problem of showing that many central values are non-zero, which can be attacked analytically if  we know two moments with a little room to spare (we gave a few examples 
of such results in the previous section), and (ii) obtaining sub-convexity bounds for $L$-functions  (there is a vast literature here, and we content ourselves to pointing to earlier surveys on this topic \cite{Fri, MicVen, Mun, IS2} and to Nelson \cite{Nelson, Nelson3} for very recent progress).

The unitary family of Dirichlet characters $(\bmod \, q)$ (for a large prime $q$) is closest in spirit to $\zeta(\frac 12 +it)$, but there are still some differences.   It is 
easy to evaluate the second moment $\sum_{\chi \pmod q}^* |L(\frac 12, \chi)|^2$ (where the $*$ indicates that the sum is restricted to primitive characters) and, in analogy with \eqref{4.1}, this is $\sim q \log q$.   The fourth moment can also be evaluated, and in analogy with Ingham's result, Heath-Brown \cite{HB5} established that  $\sum_{\chi \pmod q}^* |L(\frac 12,\chi)|^4 \sim \frac{1}{2\pi^2} q (\log q)^4$.   However an analogue  of \eqref{4.3}, obtaining lower order terms in the asymptotic formula with a ``power saving'' in the  error term, proved substantially more difficult, and was first achieved in the work of Young \cite{Young}.   Higher moments remain unknown, although one can make progress by averaging  over $q$ (see \S 6).   Another natural unitary family is the twists of a fixed Hecke eigenform by Dirichlet characters $(\bmod \, q)$.   The complexity of the second moment in this family is naively comparable to the  fourth moment of Dirichlet $L$-functions, but there 
are further formidable difficulties.  An extensive discussion of this problem, with variants and applications, may be found in the memoir of Blomer et al \cite{Blomeretal}.

 In the symplectic case of quadratic Dirichlet $L$-functions the first three moments $\sum_{|d|\le x} L(\frac 12, \chi_d)^k$ are known (see \cite{J, So5}, and for interesting work on the error term in the cubic moment see \cite{Young2, DW}), and the 
 asymptotics in these cases ($k=1, 2, 3$) take the shape of $x P_k(\log x)$ for a polynomial $P_k$ of degree $k(k+1)/2$.   We shall explain in the next section how 
 this ties in with the Keating--Snaith conjecture for the distribution of $\log L(\frac 12, \chi_d)$.   The techniques behind evaluating these moments also establish that a proportion at least $\frac 78$ of these values are non-zero (see \cite{So5}).  

As an example of an orthogonal family, consider the set ${\mathcal H}_k$ of Hecke eigenforms for the full modular group with large weight $k \equiv 0 \pmod 4$.  Here the moments $\sum_{f\in {\mathcal H}_k} L(\frac 12, f)^r$ may be evaluated for $r=1$, $2$,  and if an extra averaging over $K\le k\le 2K$ is taken, then in the cases $r=3$ and $4$ also (this follows from the techniques in \cite{IwSa}).  The asymptotic answers here are of the shape $| {\mathcal H}_k| P_r(\log k)$ for a polynomial $P_r$ of degree $r(r-1)/2$.  A sharp bound for the third moment (without an average in $k$) would provide a subconvexity bound $L(\frac 12, f) \ll k^{\frac 13+\epsilon}$, which is comparable in strength to the Hardy-Littlewood-Weyl subconvexity bound for $\zeta(\frac 12+it)$.   An analogous cubic moment (with such a subconvexity bound) has been studied in the case of Maass forms by Ivic \cite{Ivic}; interestingly, these cubic moments are also connected by a beautiful formula of Motohashi \cite{Moto} to the fourth moment of $\zeta(\frac 12+it)$.  Substantial  progress has been made towards obtaining  estimates for the fifth moment for modular forms (in the weight and level aspects) and in finding ``reciprocity relations'' among the fourth moments in different families; see 
\cite{BK, KY, Khan}.  

We mention one more striking example: the work of Conrey-Iwaniec \cite{CI} gives sharp estimates for the cubic moment of $L(\frac 12, f\times \chi)$ 
where $f$ runs over modular forms of level dividing $q$ (an odd square-free integer) and $\chi$ denotes the quadratic character $\pmod q$.  This gives a 
good Weyl-type subconvexity bound for such $L$-values, and an analogous calculation for Maass forms gives Weyl-type subconvexity bounds for quadratic Dirichlet $L$-functions (improving upon classical results of Burgess).   Further spectacular work in this direction may be found in Petrow and Young \cite{PY},  and Nelson \cite{Nelson2}.

\section{Conjectures for the asymptotics of moments} 

\noindent  Before discussing in detail the moments on the critical line, let us consider the moments on the line Re$(s)=\sigma >\frac 12$.   We mentioned in \S 2 that 
$\zeta(\sigma+it)$ is distributed like  the random object $\zeta(\sigma, {\mathbb X})$ defined in \eqref{1.4}.   We may therefore expect that for any $k\in {\mathbb N}$ and as $T\to \infty$  
\begin{equation} 
\label{5.1} 
\frac 1T \int_0^{T} |\zeta(\sigma+it)|^{2k} dt \sim {\mathbb E}[ |\zeta(\sigma, {\mathbb X})|^{2k} ]  = \sum_{n=1}^{\infty} \frac{d_k(n)^2}{n^{2\sigma}}, 
\end{equation} 
since $\zeta(\sigma, {\mathbb X})^k = \sum_{n=1}^{\infty} d_k(n){\mathbb X}(n)/n^{\sigma}$ with $d_k(n)$ being the $k$-divisor function (the series converges almost surely for $\sigma >\frac 12$).  When $\sigma >1$ it is clear that \eqref{5.1} holds (indeed for any real number $k$), since the values $|\zeta(\sigma+it)|$ lie in a compact subset of $(0,\infty)$ and the distributions match.   The case $\sigma=1$ is more delicate, but with a little more effort one can justify \eqref{5.1} here as well.
  Moving 
now into the critical strip, there is no known value of $\tfrac 12 < \sigma <1$ where 
the asymptotic \eqref{5.1} is known to hold for all $k \in {\mathbb N}$.  Indeed such a 
result would imply that $|\zeta(\sigma+it)| \ll |t|^{\epsilon}$, which  remains 
unknown for any $\tfrac 12 < \sigma <1$.  However, if one is willing to assume RH, then it is possible to approximate $\zeta(\sigma+it)^k$ by 
short Dirichlet polynomials, and then \eqref{5.1} follows for all real numbers $k$.

Returning to moments on the critical line, as mentioned previously, asymptotic formulae 
for $M_k(T)$ are known only in the cases $k=1$ and $2$.  But, using \eqref{5.1} as a 
guide,  we may guess  the order of magnitude of $M_k(T)$.  The series on
the  right side of \eqref{5.1} 
diverges when $\sigma =\tfrac 12$, but we might consider truncating that sum around size $T$.  
It is easy to show that for any real number $k$,  
\begin{equation} 
\label{5.2}
\sum_{n\le T} \frac{d_k(n)^2}{n}  \sim \frac{a_k}{\Gamma(k^2+1)} (\log T)^{k^2},
\end{equation} 
where 
\begin{equation} 
\label{5.3} 
a_k =  \prod_p \Big(1-\frac 1p \Big)^{k^2} \Big( \sum_{a=0}^{\infty} \frac{d_k(p^a)^2}{p^{a}}\Big). 
\end{equation}
Thus one might guess that for all positive real numbers $k$, $M_k(T) 
\sim C_k T(\log T)^{k^2}$ for some constant $C_k$.  Conrey and Ghosh suggested that 
it might be instructive to write $C_k$ as $g_k a_k /\Gamma(k^2+1)$, and expected 
that the unknown factor $g_k$ might have nice properties (for example, that $g_k$ would 
be a natural number when $k$ is a natural number).   The Hardy-Littlewood aymptotic 
for the second moment (see \eqref{4.2}) is in keeping with this conjecture, and gives $g_1=1$.   Similarly, Ingham's 
result on the fourth moment (see \eqref{4.3}) yields $g_2=2$.

Another way to guess at the order of magnitude for $M_k(T)$ arises from extrapolations of Selberg's central limit theorem.  
If $X$ is a random variable that is normally distributed with mean $\mu$ and variance $\sigma^2$, then 
for any real number $t$ we have 
\begin{align} 
\label{5.4} 
{\mathbb E}[e^{tX}] &= \frac{1}{\sqrt{2\pi} \sigma} \int_{-\infty}^{\infty} \exp\Big( tu -\frac{(u-\mu)^2}{2\sigma^2}\Big) du  \nonumber\\ 
&= e^{t\mu +t^2 \sigma^2/2} \frac{1}{\sqrt{2\pi } \sigma} \int_{-\infty}^{\infty} \exp\Big( - \frac{(u-\mu-t\sigma^2)^2}{2\sigma^2} \Big) du
= e^{t\mu +t^2 \sigma^2/2}. 
\end{align} 
Further the dominant contribution above comes from values of $X$ that are about $\mu + t\sigma^2  + O(\sigma)$.  Selberg's theorem tells us that 
$\log |\zeta(\tfrac 12+it)|$ is distributed like a Gaussian with mean $0$ and variance $\sim \frac 12 \log \log T$.  The calculation in \eqref{5.4} 
therefore suggests that 
$$
M_k(T) = \int_0^T \exp\Big(2k \log |\zeta(\tfrac 12+it)|\Big) dt =T \exp\Big( (2k)^2 \frac{\frac 12 \log \log T}{2}\Big) = T (\log T)^{k^2}.  
$$ 
 Moreover the dominant contribution to the $2k$-th moment should arise from values of $\zeta(\frac 12+it)$ of size $(\log T)^k$ and the set on which 
 such values are attained has measure about $T/(\log T)^{k^2}$.   We should clarify that Selberg's theorem is concerned with typical values of $\log |\zeta(\frac 12+it)|$, 
 which are on the scale of $\sqrt{\log \log T}$, whereas the moments $M_k(T)$ are concerned with the large deviations regime where $\log |\zeta(\frac 12+it)|$  
 is of size $k \log \log T$.  In this regime Selberg's result does not immediately apply, and indeed we should expect some deviations from the Gaussian, which are 
 reflected in the constant $C_k$ appearing in the conjecture for $M_k(T)$ (see \cite{Ra2, Faz}).   Later we shall discuss a coarse version of Selberg's theorem in this 
 large deviations regime \cite{So}, conditional on RH, which leads to good (conditional) upper bounds for $M_k(T)$.   To give an analogy, both $\omega(n)$ 
 (the number of distinct prime factors of $n$) and $\log d(n)/\log 2$ (with $d(n)$ being the divisor function) are additive functions that are distributed (if $n$ is 
 chosen uniformly in $[1,N]$) like a Poisson random variable with parameter $\log \log N$ --- this is the Erd{\H o}s--Kac theorem (noting that Poisson with large parameter approximates a Gaussian).  This suggests that both $\sum_{n\le N} 2^{\omega(n)}$ and $\sum_{n\le N} d(n)$ are on the scale of $N \log N$, but the 
 constants involved in the asymptotics are not immediate (and are different in the two cases).

Neither of the two heuristics given above makes a prediction for the constant $C_k = a_k g_k/\Gamma(k^2+1)$.   Indeed, until the nineties there was no clear conjecture as to the value of $g_k$ for any $k\neq 1$, or $2$.  Then Conrey and Ghosh \cite{CoGh1}, \cite{CoGh2}, based 
on an earlier conjecture of Balasubramanian, Conrey and Heath-Brown \cite{BCH-B}, advanced the 
conjecture that $g_3=42$.  A little later Conrey and Gonek \cite{CoGo}, based on 
conjectures on the asymptotics of  divisor correlation sums (as in \eqref{4.5}), 
arrived again at the conjecture that $g_3=42$ (see Ng \cite{Ng} for further  work on making this precise), while also advancing the conjecture that $g_4=24024$.  These methods did not extend to produce conjectures for larger $k$, and the problem once again seemed stuck.  A great advance was made when Keating and Snaith \cite{KeSn2}, 
using ideas from random matrix theory, arrived at the following remarkable conjecture 
for $M_k(T)$ for all positive real numbers $k$.

\begin{conjecture}[Keating and Snaith]  For any positive real number $k$, as $T\to \infty$, we have 
$M_k(T) \sim g_k \frac{a_k}{\Gamma(k^2+1)} T (\log T)^{k^2}$, 
with
$$
g_k = \Gamma(k^2+1) \frac{G(1+k)^2}{G(1+2k)}, 
$$ 
where $G$ is the Barnes $G$-function.   In particular, if $k \in {\mathbb N}$ then 
$$
g_k = (k^2)! \prod_{j=0}^{k-1} \frac{j!}{(k+j)!}, 
$$  
so that $g_1=1$, $g_2=2$, $g_3=42$, and $g_4=24024$. 
\end{conjecture}

We recall that the Barnes $G$-function is an entire function of order $2$ which 
satisfies the functional equation $G(z+1) = \Gamma(z) G(z)$ 
with the normalization $G(1)=1$.  Thus for a natural number $n$, one has $G(n) 
= \prod_{j=0}^{n-2} j!$.   

The key insight of Keating and Snaith was to quantify and develop in the context of value distribution problems a conjectural connection between 
the distribution of zeros of the Riemann zeta function and the distribution of eigenvalues of large random matrices.   Nearly fifty years back, Montgomery \cite{Mo} 
initiated a study of the spacings between the ordinates of zeros of the Riemann zeta function, and a chance conversation with Dyson revealed that his partial results 
on this question matched corresponding statistics in the study of spacings between eigenvalues of large random matrices.  Assuming RH for clarity, let $\gamma_1 \le \gamma_2 \le \ldots$ denote the sequence of non-negative ordinates of zeros of $\zeta(s)$ (written with multiplicity), so that from \eqref{2.25} it follows that $\gamma_n \sim 2\pi n/\log n$.  The question then is to determine the distribution (as $n\to \infty$) of $(\gamma_{n+1}-\gamma_n) (\log \gamma_n)/(2\pi)$, which has been normalized to have mean spacing $1$.  For example, with what frequency does this normalized spacing lie in a given interval $(\alpha, \beta) \subset (0,\infty)$?  One way to express the (amazing!) conjectured answer is as follows.  Consider a random element $g$ drawn from the unitary group $U(N)$ with respect to the Haar measure $dg$ (normalized so that $U(N)$ has volume $1$).  Each such $g$ has eigenvalues $e^{i\theta_1}$, $\ldots$, $e^{i\theta_N}$ with the angles ordered $0 \le \theta_1 \le \theta_2\le \ldots \le \theta_N <2\pi$, and consider the spacings $(\theta_{n+1}-\theta_n) N/(2\pi)$ (normalized to have average approximately $1$).  Average this spacing distribution over the whole group $U(N)$, and finally let $N\to \infty$.  For example, we could count the frequency with which $(\theta_{n+1}-\theta_n)N/(2\pi)$ lies in $(\alpha,\beta)$, average that frequency over $U(N)$, and take the limiting frequency as $N\to \infty$.    
The model that we have described is known as the \emph{Circular Unitary Ensemble} (CUE), and the same distribution for nearest neighbor spacings arises in 
other models of random matrices such as the Gaussian Unitary Ensemble (GUE).  

Theoretical support for this link between zeros of $\zeta(s)$ and random matrix theory arose first with Montgomery's calculation of the \emph{pair correlation} of 
zeros (in certain ranges) mentioned earlier, and this was generalized to general $n$-level correlations in the work of Rudnick and Sarnak \cite{RuSa}.  Experimental support for this link comes from extensive computations of Odlyzko \cite{Odl} who considered the spacing distribution of about $175$ million zeros around the $10^{20}$-th zero (which occurs at height around $1.5 \times 10^{19}$), and found an astonishingly close match between the empirical data and the 
predicted answer.  Yet, Odlyzko's data found that the numerical data did not match closely some other statistics for $\zeta(s)$ such as Selberg's theorem on $\log \zeta(\frac 12+it)$. One might attribute such deviations to the slow growth of the variance $\log \log T$, which even at height $10^{19}$ is only about $3.7$, but Keating and Snaith \cite{KeSn2} suggested a much more insightful explanation.   They posited that properties of $\zeta(\frac 12+it)$ for $t$ around a specific height $T$ may be modeled by analogous objects for random matrices of a specific size $N$, refining the expectation that the large $T$ and large $N$ limits coincide.   The relation between $N$ and $T$ is suggested by the average spacing between the zeros at height $T$, which is about $(2\pi)/\log (T/2\pi)$ by \eqref{2.25}, and the average spacing between 
eigenvalues, which is about $(2\pi)/N$.  Setting these equal, we arrive at the correspondence $N \approx \log (T/2\pi)$.   The analogue of the zeta function, which is determined by its zeros, is the characteristic polynomial of a random matrix, which is determined in a similar fashion by its eigenvalues.  Keating and Snaith determined the distribution of $\log \text{det} (I - ge^{-i\theta})$ for a random matrix $g\in U(N)$, and found that in the large $N$ limit this tends to a complex Gaussian 
with mean $0$ and variance $\log N$ (analogously to Selberg's theorem), but there are lower order terms that are significant for finite $N$.   The range of Odlyzko's computations, $T \approx 1.5 \times 10^{19}$, corresponds to matrices of size $N =42$, and Keating and Snaith found an excellent fit between Odlyzko's numerical data for $\log \zeta(\frac 12+it)$ and the distribution of $\log \text{det} (I - ge^{-i\theta})$ for random $g\in U(42)$ (see Figure 1 in \cite{KeSn2}).

Returning to the moments, one might now hope to understand the asymptotic behavior of $M_k(T)$ by computing the analogous moments in the context of $U(N)$: namely 
\begin{equation} 
\label{5.5} 
\int_{g\in U(N)} \frac{1}{2\pi } \int_0^{2\pi} | \text{det} (I- ge^{-i\theta})|^{2k} d\theta dg = \int_{g\in U(N)} |\text{det}(I-g)|^{2k} dg. 
\end{equation} 
By the Weyl integration formula expressing the measure $dg$ in terms of the eigenvalues of $g$, this equals the multiple integral 
\begin{equation} 
\label{5.6} 
\frac{1}{(2\pi)^N N!} \int_0^{2\pi} \cdots \int_0^{2\pi} \Big|\prod_{j=1}^{N} (1-e^{i\theta_j})\Big|^{2k} \prod_{1\le j< m \le N} |e^{i\theta_j}-e^{i\theta_m}|^2 
d\theta_1 \cdots d\theta_N. 
\end{equation} 
It turns out that the integral in \eqref{5.6} may be evaluated exactly using a remarkable formula of Selberg \cite{Sel3} (see \cite{FW} for many developments arising from the Selberg integral) and it equals 
\begin{equation} 
\label{5.7} 
\prod_{j=1}^{N} \frac{\Gamma(j) \Gamma(2k+j)}{(\Gamma(j+k))^2}  \sim g_k \frac{N^{k^2}}{\Gamma(k^2+1)}, 
\end{equation} 
where $g_k$ is as in Conjecture 5.1, and the asymptotic holds for large $N$.   The constant $g_k$ has an intriguing combinatorial interpretation as 
the number of standard Young tableaux of shape $k\times k$ (that is, the number of ways of filling a $k\times k$ array with the numbers $1$, $\ldots$, $k^2$ 
such that the entries along each row and column are in increasing order).  See \cite{BuGa, DeH, KRRR} for related combinatorial discussions, and 
\cite{CF, Ger} for discussions on the divisibility properties of $g_k$ and related constants.   

This calculation motivates Conjecture 5.1, but note that no primes appear in the random matrix model, and so the constant $a_k$ must be ``put in by hand.''  Here we note that the Euler product for $a_k$ in \eqref{5.3} arises naturally upon considering $\E [ |(1-X(p)/\sqrt{p})|^{-2k}] = \sum_{a=0}^{\infty} d_k(p^a)^2/p^{a}$ with $X(p)$ 
chosen uniformly from the unit circle.  Thus the constant $a_k$ may be thought of as arising from a version of the random Euler product, while the $g_k$ term arises from the local behavior of zeros of the zeta function.  For an exploration of Conjecture 5.1 along these lines, developing a hybrid Euler--Hadamard product, see 
the work of Gonek, Hughes, and Keating \cite{GHK}.  We mentioned earlier the analogy with determining asymptotics for multiplicative functions such as 
$k^{\omega(n)}$ or $d_k(n)$, and here the known asymptotic formulae (going back to Landau, Selberg, and Delange) factor as a ``local'' product over primes together with a ``global'' term determined by the Poisson behavior of $\omega(n)$; for an interesting discussion of this analogy see \cite{JKN}.

Random matrix theory also informs our understanding of moments of central values of $L$-functions in families.  While the distribution of spacings between zeros at large height for any given $L$-function is  expected to follow the same law that we described for $\zeta(s)$ (see \cite{RuSa}), the distribution of the zeros close to the central point $\frac 12$ can vary depending on the particular family.   Based on analogies with the function field case, Katz and Sarnak \cite{KaSa} found (conjecturally) that the distribution of zeros near $\frac 12$ in families of $L$-functions fell into the three categories \emph{unitary, symplectic, and orthogonal} (which we discussed in 
\S 3), and that these distributions matched the distribution of the eigenvalues close to $1$ of large random matrices chosen from $U(N)$, $USp(2N)$, or $SO(2N)$ 
(or $SO(2N+1)$  depending on the sign of the functional equation).  To give an illustration of the Katz--Sarnak conjectures, consider the family of quadratic Dirichlet $L$-functions $L(s, \chi_d)$ as $d$ ranges over fundamental discriminants, which is expected to have symplectic symmetry.   The density of zeros of $L(s,\chi_d)$ near $\frac 12$ is about $(\log |d|)/(2\pi)$, and a sample question is to understand the distribution of $\gamma_1 \frac{\log |d|}{2\pi}$ where $\gamma_1$ is 
the least non-negative ordinate of a zero of $L(s,\chi_d)$.  To describe the conjectured answer, consider a random matrix $g \in USp(2N)$ (chosen with respect to Haar measure normalized to have total volume $1$) and write its eigenvalues as $e^{\pm i\theta_1}$, $e^{\pm i\theta_2}$, $\ldots$, $e^{\pm i \theta_N}$ with $0\le \theta_1 \le \ldots \le \theta_N \le \pi$.  Then as $d$ varies over fundamental discriminants $|d|\le X$ with $X \to \infty$, the distribution of $\gamma_1 \frac{\log |d|}{2\pi}$ 
is identical to the limiting distribution of $\theta_1 \frac{2N}{2\pi}$ for randomly chosen $g\in USp(2N)$ as $N\to \infty$.  

Conrey and Farmer \cite{CF} proposed that the moments of central values of $L$-functions in families are also dictated by the symmetry type in the Katz--Sarnak conjectures.   In particular, the analogue of the factor $g_k$ should depend only on the symmetry type and not on the particular family, whereas the analogue of the factor $a_k$ will be sensitive to the particular family (in a straightforward way).   This was developed further by Keating and Snaith \cite{KeSn1}, who modeled properties of the central $L$-values by the characteristic polynomial $\text{det}(I-ge^{-i\theta})$ evaluated at $\theta=0$, with the size parameter $N$ of the random matrix ensemble chosen to match with the density of zeros in the family.  Indeed it is a consideration of the behavior of $\log \text{det}(I-g)$ in $USp(2N)$ or $SO(2N)$ that informed their conjectures for the analogues of Selberg's theorem in symplectic and orthogonal families (discussed in \S 3).  

Just as extrapolating Selberg's theorem allows us to guess the order of 
magnitude of moments of $\zeta(s)$, the Keating--Snaith log normality conjectures together with the calculation in \eqref{5.4} gives an understanding of the order of magnitude of moments in families.  For example, in the symplectic example of moments of $L(\frac 12,\chi_d)$ with $|d|\le X$, since $\log L(\frac 12,\chi_d)$ is 
conjectured to be normal with mean $\sim \frac 12 \log \log X$ and variance $\sim \log \log X$, the calculation in \eqref{5.4} suggests that $\sum_{|d| \le X} L(\frac 12,\chi_d)^k$ is of size $X (\log X)^{\frac{k(k+1)}{2}}$.   Similarly in the orthogonal case of Hecke eigenforms $f\in {\mathcal H}_k$, since $\log L(\frac 12,f)$ is 
expected to be normal with mean $\sim -\frac 12\log \log k$ and variance $\sim \log \log k$, the moments $\sum_{f \in {\mathcal H}_k} L(\frac 12, f)^r$ may be expected to be of order $k (\log k)^{\frac{r(r-1)}{2}}$.  

Further, by considering moments of $\text{det}(I-g)$ in the appropriate matrix group, Keating and Snaith \cite{KeSn1} formulated analogues of Conjecture 5.1 
in families of $L$-functions.   For example, in the family of quadratic Dirichlet $L$-functions $L(s,\chi_d)$, the analogue of the constant $g_k$ is predicted by 
considering 
\begin{align*}
\int_{g\in USp(2N)} \text{det}(I-g)^{k} dg = 2^{2Nk} \prod_{j=1}^{N} \frac{\Gamma(1+N+j) \Gamma(1/2+ k+j)}{\Gamma(1/2+j) \Gamma(1+k+N+j)} 
\sim f_k \frac{N^{k(k+1)/2}}{\Gamma(k(k+1)/2+1)}.
\end{align*}
This calculation again reduces to the Selberg integral, and the constant $f_k$ may be expressed in terms of the Barnes $G$-function.  If $k$ is a natural number then $f_k$ takes the pleasant form $(k(k+1)/2)!/ \prod_{j=1}^{k} (2j-1)!!$.  After incorporating an analogue of the constant $a_k$ in \eqref{5.3}, which here is (with ${\mathbb X}(p)$ denoting the random variables modeling quadratic characters discussed in \S 1) 
\begin{align*}
&\prod_p \Big(1-\frac 1p \Big)^{\frac{k(k+1)}{2}} \E \Big[ \Big(1- \frac{{\mathbb X}(p)}{\sqrt{p}}\Big)^{-k}\Big] \\
= &\prod_{p}  \Big(1-\frac 1p \Big)^{\frac{k(k+1)}{2}} \Big( \frac{p}{2(p+1)} \Big(\Big( 1+\frac{1}{\sqrt{p}}\Big)^{-k} + \Big( 1- \frac 1{\sqrt{p}}\Big)^{-k} \Big) + \frac{1}{p+1}\Big),
\end{align*}
we arrive at a conjecture for the moments of $L(\frac 12,\chi_d)$, which matches  the known asymptotics for the first three moments.

The Keating--Snaith conjectures identify the leading order term in the asymptotics for moments, but there will be lower order terms (just a logarithm smaller) which are not identified.   We may see this already in the asymptotics for the second and fourth moments of $\zeta(\frac 12+it)$ (see \eqref{4.2} and \eqref{4.3}), and other examples in families given in \S 4.   Identifying such lower order terms is of interest because the leading order constant in Conjecture 5.1,  $a_kg_k/\Gamma(k^2+1)$ tends rapidly to zero as $k$ grows, and so for the ranges of $T$ in which numerical investigations may be carried out, the lower order terms may dominate the eventual main term.   When $k$ is a positive integer, Conrey et al \cite{CFKRS} conjectured that $M_k(T) = \int_0^T P_k(\log t/2\pi) dt +O(T^{1-\delta})$ (for some $\delta>0$, and perhaps even any $\delta< \frac 12$ is permissible) for a polynomial $P_k$ of degree $k^2$ with leading coefficient $a_k g_k/(k^2!)$, and they gave  a ``recipe'' for determining all the coefficients of $P_k$.  Their recipe predicts the full main term for integral moments in many families of $L$-functions, but it remains open to give an asymptotic expansion when $k$ is not an integer.  The paper \cite{CFKRS} also gives numerical evidence towards the full moment conjecture, and further data may be found in \cite{HO}.    A related approach via \emph{multiple Dirichlet series} is described in the work of Diaconu, Goldfeld, and Hoffstein \cite{DGH} who develop conjectures for the integral moments of quadratic Dirichlet $L$-functions (which are in agreement with \cite{CFKRS}).  

We give a brief illustration of the recipe from \cite{CFKRS} in the unitary family of Dirichlet $L$-functions $\chi\pmod q$ with $q$ a large prime.  For simplicity, we consider only even characters (thus $\chi(-1)=1$), where the functional equation reads $\Lambda(s,\chi)= (q/\pi)^{s/2} \Gamma(s/2) L(s, \chi) = \epsilon_{\chi} \Lambda(1-s,\overline{\chi})$ with $\epsilon_{\chi}$ satisfying $|\epsilon_{\chi}|=1$ and $\epsilon_{\chi} \epsilon_{\overline{\chi}} =1$.  Let $\underline{\alpha} = (\alpha_1,\ldots, \alpha_k)$, and $\underline{\beta}  = (\beta_1,\ldots, \beta_k)$ denote two $k$-tuples of complex numbers (thought of as small), and we also find it convenient to write $\alpha_{k+j} =\beta_j$ and think of $(\underline{\alpha},\underline{\beta})$ as the $2k$-tuple $(\alpha_1, \ldots, \alpha_{2k})$.  Instead of considering $|L(\frac 12, \chi)|^{2k}$ directly, we work with 
$$
\Lambda(\chi;\underline{\alpha},\underline{\beta}) := \prod_{j=1}^{k} \Lambda(\tfrac 12+\alpha_j, \chi) \Lambda(\tfrac 12 -\beta_j, \overline{\chi})
$$
 and finally let all the parameters $\alpha_j$ and $\beta_j$ tend to zero (which would then equal $|L(\frac 12,\chi)|^{2k}$ multiplied by the constant $(q/\pi)^{k/2} \Gamma(1/4)^{2k}$).   Permuting the $k$ entries in $\underline{\alpha}$, 
or the $k$ entries in $\underline{\beta}$ does not change $\Lambda(\chi;\underline{\alpha}, \underline{\beta})$.  Less obviously, it turns out that $\Lambda(\chi; \underline{\alpha},\underline{\beta})$ is invariant under any permutation of the $2k$-entries in $(\underline{\alpha}, \underline{\beta})$; this is because any such permutation must change some $\ell$ of the $\alpha$'s to $\beta$'s and a corresponding number of $\beta$'s to $\alpha$'s and $2\ell$ applications of the 
functional equation ($\ell$ of them with a factor $\epsilon_\chi$ and $\ell$ with a factor $\epsilon_{\overline{\chi}}$) justify the claim.   Thus any conjecture that 
we propose for $\sum_{\chi} \Lambda(\chi;\underline{\alpha},\underline{\beta})$ must satisfy this $S_{2k}$ symmetry. 

Now if Re$(s)$ is large, expanding the $L$-functions into their Dirichlet series, we may write 
\begin{align}
\label{5.8}
\prod_{j=1}^{k} \Lambda(s+\alpha_j,\chi) \Lambda(s-\beta_j,\chi) &= 
\prod_{j=1}^{k} \Big(\frac{q}{\pi}\Big)^{s+\frac{\alpha_j-\beta_j}{2}} \Gamma\Big(\frac{s+\alpha_j}{2}\Big)\Gamma\Big(\frac{s-\beta_j}{2}\Big) \nonumber \\
&\times 
\sum_{m, n=1}^{\infty} \frac{\sigma(m;\underline{\alpha})}{m^s} \chi(m)  \frac{\sigma(n;-\underline{\beta})}{n^s} \overline{\chi}(n), 
\end{align}
where $\sigma(m;\underline{\alpha}) = \sum_{m=m_1 \cdots m_k} m_1^{-\alpha_1}\cdots m_k^{-\alpha_k}$ and similarly $\sigma(n;-\underline{\beta}) 
= \sum_{n=n_1\cdots n_k} n_1^{\beta_1} \cdots n_k^{\beta_k}$, so that if $\alpha_i=\beta_i=0$ these would simply be the $k$-divisor function.  We average this over all the even characters $\bmod \, q$ (omitting the trivial character), and hypothesize that only the diagonal terms $m=n$ survive this averaging.  This is of course not justified, but is similar to the first heuristic we gave in this section for the order of magnitude of moments.  After a computation with Euler products, these terms give (for the sum over $m, n$ in \eqref{5.8}) 
\begin{equation} 
\label{5,9}
\sum_{n=1}^{\infty} \frac{\sigma(n;\underline{\alpha})\sigma(n;-\underline{\beta})}{n^{2s}} 
= {\mathcal A}(s;\underline{\alpha},\underline{\beta}) \prod_{j, \ell =1}^{k} \zeta(2s+\alpha_j-\beta_\ell),
\end{equation} 
where ${\mathcal A}$ is given by an Euler product that converges absolutely in Re$(s)> \frac 12 -\delta$ if $\alpha_j$, $\beta_j$ are small enough.  This factor ${\mathcal A}$ is similar to the $a_k$ appearing in \eqref{5.3}.  Evaluating this at $s=\frac12$, we arrive at a candidate for the average value of $\Lambda(\chi;\underline{\alpha},\underline{\beta})$, namely 
\begin{equation} 
\label{5.10} 
{\mathcal C}(\underline{\alpha}, \underline{\beta}) = \prod_{j=1}^{k} \Big(\frac{q}{\pi}\Big)^{\frac{1+\alpha_j-\beta_j}{2}} \Gamma\Big(\frac{\frac 12+\alpha_j}{2}\Big)\Gamma\Big(\frac{\frac 12-\beta_j}{2}\Big) 
{\mathcal A}(\tfrac 12;\underline{\alpha},\underline{\beta}) \prod_{j, \ell =1}^{k} \zeta(1+\alpha_j-\beta_\ell).
\end{equation}

The candidate answer ${\mathcal C}(\underline{\alpha},\underline{\beta})$ is invariant when the entries of $\underline{\alpha}$ are permuted, or 
when the entries of $\underline{\beta}$ are permuted, but does not have the $S_{2k}$ symmetry we require of being allowed to permute the $2k$-entries of 
$(\underline{\alpha},\underline{\beta})$.   The beautifully simple answer proposed in \cite{CFKRS} is to symmetrize ${\mathcal C}(\underline{\alpha},\underline{\beta})$ 
by summing over all $\binom{2k}{k}$ cosets of $S_{2k}/(S_k \times S_k)$: 
\begin{equation} 
\label{5.11} 
\sum_{\pi \in S_{2k}/S_k\times S_k} {\mathcal C}(\pi(\underline{\alpha}, \underline{\beta})). 
\end{equation} 
While the expression in \eqref{5.10} has singularities whenever $\alpha_j =\beta_\ell$, the symmetrized expression in \eqref{5.11} turns out to be 
regular when $|\alpha_j|, |\beta_j|$ are small.  Now setting $\alpha_1 = \cdots = \alpha_k = \beta_1 =\cdots = \beta_k =0$, we arrive at the conjectured answer for the 
average of $|L(\frac 12,\chi)|^{2k}$.  The leading term matches the Keating--Snaith conjecture, but now we also have the full polynomial of degree $k^2$.

To end our discussion of the moment conjectures, we mention recent work of Conrey and Keating \cite{CKV} which aims to give a heuristic derivation of 
the moment conjectures of $\zeta(s)$ from correlations of divisor functions (as in \cite{CoGo} for the sixth and eighth moments).  It would be of interest to develop their work in other families of $L$-functions.   Sawin \cite{Saw} develops a heuristic approach based on representation theory which (conditional on some 
hypotheses) recovers the recipe in Conrey \emph{et al} \cite{CFKRS} in the function field setting (with a fixed field of constants).

\section{Progress towards understanding the moments} 

\noindent In \S 4 we gave a number of examples where asymptotics for low moments are known, and all of these are in agreement with the conjectures described in the previous section.  A rule of thumb suggests that an asymptotic for a moment may be computed if there are more elements in the family compared to the complexity of approximating the required power of the $L$-value (what we have informally called the complexity can be thought of as the square-root of the  \emph{analytic conductor}, see \cite{IS2}).  For example, as we saw in \eqref{4.4} $\zeta(\frac 12+it)$ may be approximated by (two) Dirichlet polynomials of length about $\sqrt{t}$, allowing for the calculation of the second and fourth moments.   This rule of thumb is only a rough guide, and can be difficult to attain.   For example, the fourth moment of Dirichlet $L$-functions $\bmod\, q$ (evaluated in \cite{Young}), or the mean square of  twists of a modular form by Dirichlet characters $\bmod\, q$ (see \cite{Blomeretal, KMS}) may seem of comparable difficulty to the fourth moment of the zeta function, but the first two problems turn out to be substantially harder.  The largest moment that may be computed by this rule of thumb recovers the convexity bound for the $L$-value, and so there is great interest in going beyond this range, either by shrinking suitably the family over which we average, or by adding an extra short Dirichlet polynomial to the moment.  

From the viewpoint of verifying the moment conjectures (for example to check the constants $42$ and $24024$ appearing in the sixth and eighth moments) one might look for large families where the complexity is still small.    The family of primitive Dirichlet characters $\chi \pmod q$ ranging over all moduli $q\le Q$ is a good example, where the size of the family is about $Q^2$ whereas the complexity of such $L(\frac 12, \chi)$ is about $\sqrt{Q}$.  This suggests the possibility of evaluating the sixth and eighth moments in this family, and indeed the large sieve gives a quick upper bound of the correct order of magnitude for these moments (see 
\cite{Hux}).  By developing an asymptotic version of the large sieve, Conrey, Iwaniec and Soundararajan \cite{CIS} obtained an asymptotic formula for 
\begin{equation} 
\label{6.1}
\sum_{q\le Q}\, \,   \sum_{\chi (\bmod \, q)}^{\flat} \int_{-\infty}^{\infty} |\Lambda(\tfrac 12+iy,\chi)|^6 dy, 
\end{equation} 
where $\Lambda(s,\chi)=(q/\pi)^{s/2} \Gamma(s/2)L(s,\chi)$ denotes the completed $L$-function, and the $\flat$ indicates a sum over even primitive characters $\chi$.   Here the averaging over $y$ is a technical defect, needed for the proof,  which (owing to the rapid decay of the $\Gamma$-function) may be thought of as an integral over essentially a bounded range of $y$.  This asymptotic formula verified the predicted constant $g_3=42$ in this instance, and moreover \cite{CIS} obtained a similar asymptotic formula with shifts $(\alpha_1,\alpha_2,\alpha_3)$ and $(\beta_1,\beta_2,\beta_3)$ which verified the recipe of \cite{CFKRS} in this situation and yielded the full polynomial of degree $9$ in $\log Q$ for \eqref{6.1}.  Chandee and Li \cite{ChLi1} tackle the analogue of \eqref{6.1} for the eighth moment, and obtain an asymptotic 
formula conditional on the Generalized Riemann Hypothesis.   Their work confirmed that $g_4=24024$ in this instance, but they could only verify the leading order term in the asymptotic and not the full polynomial of degree $16$.  Forthcoming work of Chandee, Li, Matom{\" a}ki and Radziwi{\l \l} (see \cite{Ra3} for an announcement) removes the imperfection of the average over $y$ in \eqref{6.1} for the sixth moment while still obtaining the full asymptotic formula with power saving.  They also obtain the leading order behavior of the eighth moment without invoking GRH, and without the integral in $y$.  

The family of  newforms of a fixed weight $k$ for the group $\Gamma_1(q)$ with $q$ a large prime  offers another instance of a large family where the complexity (or analytic conductor) remains small.  These correspond to newforms for $\Gamma_0(q)$ with character $\chi\pmod q$.  This is a family of about $q^2$ elements, and is unitary since almost all of the characters $\chi \pmod q$ are not real.   The complexity of the $L$-values is about size $\sqrt{q}$, and we may hope to address the sixth and eighth moments.   Chandee and Li \cite{ChLi2} give an asymptotic for the sixth moment analogous to \eqref{6.1} in this family (confirming again $g_3=42$), and obtain in \cite{ChLi3} a good upper bound for the eighth moment.  It would be of interest to find further examples of families where one can compute higher moments, and in particular to obtain such examples of symplectic and orthogonal families.  The recent work of Nelson 
\cite{Nelson3} on subconvexity for automorphic $L$-functions raises the hope that one might be able to compute high moments in $GL(n)$ families 
for suitably large $n$.

\smallskip 

In addition to examples where asymptotics for moments are known, substantial progress has been made in obtaining upper and lower bounds of the conjectured 
order of magnitude in a good deal of generality.   Summarizing the work of many researchers, here is our knowledge of such bounds for the moments of $\zeta(\frac 12+it)$. 

\begin{theorem}  Let $k >0$ and $T\ge e$  be real numbers.  Then there are positive constants $c_k$ and $C_k$ such that 
$$ 
c_k T(\log T)^{k^2} \le \int_0^{T} |\zeta(\tfrac 12+it)|^{2k} dt \le C_k T(\log T)^{k^2}. 
$$ 
Here the lower bound holds unconditionally for all $k$, while the upper bound holds unconditionally in the range $0 < k\le 2$, and the upper bound holds assuming the truth of the Riemann Hypothesis for all $k >2$.  
\end{theorem} 

We shall now discuss this result and its extensions in families of $L$-functions.  The discussion splits naturally into three parts (i) lower bounds for moments, (ii) unconditional upper bounds for moments, and (iii) upper bounds assuming RH or GRH.  

The lower bound stated in Theorem 6.1 was first established by Ramachandra \cite{Ram1, Ram2} in the case when $2k$ is a natural number.  This was 
then extended by Heath-Brown \cite{HB3} to the case when $k$ is any positive rational number, but the constants $c_k$ in his result depended upon the 
height of the rational number $k$, so that the method did not extend to irrational $k$.   Further, the techniques in these works were specific to the ``$t$-aspect'' and  did not extend to moments in families of $L$-functions.  Rudnick and Soundararajan \cite{RuSo1, RuSo2} developed an alternative approach, which worked in general families.  For example, their method would show that $\sum_{|d|\le X} |L(\frac 12, \chi_d)|^{k} \ge c_k X(\log X)^{k(k+1)/2}$ for all rational $k\ge 1$ and a suitable positive constant $c_k$, which again did not vary continuously with $k$ but depended on the height of the rational number $k$.  This was further refined by Radziwi{\l \l} and Soundararajan \cite{RaSo}, who obtained the lower bounds in Theorem 6.1 for all real $k \ge 1$ with $e^{-30k^4}$ being a permissible value for $c_k$ if $T$ is large. A further round of simplification is carried out in Heap and Soundararajan \cite{HeapSo}, which also gives the lower bound in Theorem 6.1 for real $0< k\le 1$.

The story for lower bounds may be encapsulated in the following broad principle.   Whenever we can compute the mean value of $L(\frac 12)$ multiplied by short Dirichlet polynomials in a family, we can obtain lower bounds of the right order of magnitude for the moments $|L(\frac 12)|^k$ for any real $k \ge 1$.  Of course, in general H{\" o}lder's inequality will give lower bounds for higher moments in terms of smaller moments, but those would not be of the conjectured order of magnitude since the exponent of the logarithm in the moment conjectures is quadratic in $k$.  If we can also  compute the mean value of $|L(\frac 12)|^2$ multiplied by short Dirichlet polynomials, then we can obtain lower bounds of the right order of magnitude for the 
moments $|L(\frac 12)|^k$ in the range $0< k \le 1$ as well.  It may seem puzzling why the problem for small $k$ should require more information than for large $k$, but in fact this is natural.  Consider letting $k \to 0^+$.  Then the moments $|L(\tfrac 12)|^k$ essentially pick up whether $L(\tfrac 12)$ is zero or not, so that lower bounds for the small moments encode lower bounds for non-vanishing.  The analytic methods for producing non-zero values of $L(\frac 12)$ (the \emph{mollifier method}) rely on knowledge of the first two moments in the family (with a little room to spare).  Thus we may establish (using the methods of either \cite{RaSo} or \cite{HeapSo}) that for all real $k>0$, 
\begin{equation} 
\label{6.2} 
\sum_{\chi \pmod q}|L(\tfrac 12, \chi)|^{2k} \gg_k q (\log q)^{k^2},
\end{equation} 
where $q$ is a large prime, and that 
\begin{equation} 
\label{6.3} 
 \sum_{|d|\le X} |L(\tfrac 12,\chi)|^k \gg_k X (\log X)^{\frac{ k(k+1)}{2}}.   
\end{equation} 
In the family of quadratic twists of a fixed Hecke eigenform $f$, we only have access to the first moment and not the second, and therefore we only know in the 
range $k\ge 1$ that 
\begin{equation} 
\label{6.4} 
\sum_{|d|\le X} L(\tfrac 12, f\times \chi_d)^k \gg_k X(\log X)^{\frac{k(k-1)}{2}}. 
\end{equation}

We now turn to the unconditional upper bounds in Theorem 6.1, which were established in the special cases $k=1/n$ or 
$k=1+1/n$ (for natural numbers $n$) by Heath-Brown \cite{HB3} and Bettin, Chandee, and Radziwi{\l \l} \cite{BCR} respectively.  Then in Heap, Radziwi{\l \l}, and Soundararajan \cite{HRS} the bound was established for all $0< k\le 2$, as an illustration of an upper bound principle (complementing the one for lower bounds above) 
enunciated by Radziwi{\l \l} and Soundararajan \cite{RaSo2}.  Whenever we can compute a moment $|L(\frac 12)|^k$ (usually with $k$ being a positive integer) 
together with flexibility to introduce a short Dirichlet polynomial, we can obtain upper bounds of the conjectured order of magnitude for the moments $|L(\frac 12)|^r$ 
for all $0 < r\le k$.  Thus one can obtain complementary upper bounds in \eqref{6.2} for $k\le 1$ (with more effort, using Young's work \cite{Young}, this 
could perhaps be extended to the range $k\le 2$), matching upper bounds in \eqref{6.3} in the range $k\le 2$ (if one knew the positivity of $L(\frac 12,\chi_d)$ this would also follow in the 
range $k\le 3$ and it would be interesting to attain that range unconditionally), and for the family in \eqref{6.3} for $k\le 1$ (this is the example carried out in \cite{RaSo2}).

 The conditional bounds in Theorem 6.1 originated from work of Soundararajan \cite{So} who established (assuming RH)  the nearly sharp bound $M_k(T) \ll_{k,\epsilon} T(\log T)^{k^2+\epsilon}$.  This was then refined in the beautiful work of Harper \cite{Harper3} to its present sharp form.   The method is very general and applies in any family where we are able to compute the mean values of short Dirichlet polynomials.   Thus (assuming GRH in the relevant families) 
 one can obtain upper bounds of the correct order of magnitude for all non-negative $k$ in the examples \eqref{6.2}, \eqref{6.3}, and \eqref{6.4}.
 
 The main idea behind the conditional upper bounds in Theorem 6.1 is that on RH (or GRH) one can obtain an \emph{upper bound} for 
 $\log |\zeta(\tfrac 12+it)|$ (or more  generally the logarithm of central $L$-values) just in terms of sums over primes.  This is related to the ideas behind 
 Selberg's central limit theorem and the one sided versions for $L$-values that we discussed in Sections 2 and 3.  A barrier to approximating $\log |\zeta(\frac 12+it)|$ by a suitable Dirichlet polynomial is the presence of zeros near $\frac 12+it$; the crucial point is that these zeros should only make $|\zeta(\frac 12+it)|$ smaller, so that 
such Dirichlet polynomials could serve as an upper bound.   One way to see this is to note that RH is equivalent to the property that, with $s=\sigma+it$ 
$$ 
| \xi(s)| = \Big| s(s-1)\pi^{-s/2} \Gamma(s/2)\zeta(s) \Big| = \prod_{\rho} \Big| 1-\frac{s}{\rho}\Big| 
$$ 
is an increasing function of $\sigma$ in $\sigma \ge \frac 12$ for any fixed $t$.   This permits bounding $|\zeta(\frac 12+it)|$ in terms of $|\zeta(\sigma_0 +it)|$ for 
any $\sigma_0 > \frac 12$, and one can adapt Selberg's ideas to approximate $\log |\zeta(\sigma_0+it)|$.  In this manner, it was shown in \cite{So} that 
for $T\le t\le 2T$ and any $2\le x\le T^2$ one has, assuming RH and with $\sigma_0= \frac 12 + \frac{1}{\log x}$, 
\begin{equation} 
\label{6.5} 
\log |\zeta(\tfrac 12+it)| \le \text{Re }\sum_{2\le n\le x} \frac{\Lambda(n)}{n^{\sigma_0 +it}\log n} \frac{\log (x/n)}{\log x} + \frac{\log T}{\log x} + O\Big( \frac{1}{\log x}\Big). 
\end{equation} 
Analogous bounds hold for $\log |L(\tfrac 12)|$ if a corresponding GRH is assumed.   

The usefulness of \eqref{6.5} lies in its flexibility with choosing the parameter $x$.  If $x$ is suitably small, then the distribution of the 
sum (which is essentially Re$\sum_{p\le x} 1/p^{\frac 12+it}$)  in \eqref{6.5} can be understood accurately by studying its moments (as we discussed in \S 2 and \S 3), but we lose some information in the $\log T/\log x$ term.  Here it is also useful to split the sum over $p$ into different ranges (say $p\le z$ and $z<p\le x$); for small ranges of $p$, more moments may be computed so that a finer understanding of the sum is possible, while for 
the larger ranges the slow growth of the variance (which is roughly $\sum_{z\le p \le x} 1/p \sim \log \frac{\log x}{\log z}$) permits a good understanding with fewer moments.  In this way \cite{So} established a coarse version of Selberg's central limit theorem in the large deviations regime, showing that in the 
range $\sqrt{\log \log T} \le V = o(\log \log T \log \log \log T)$ one has 
\begin{equation} 
\label{6.6} 
\text{meas}\Big\{ T\le t\le 2T: \log |\zeta(\tfrac 12+it)| \ge V \Big\} \ll T \exp\Big( - \frac{V^2}{\log \log T} (1+o(1))\Big). 
\end{equation} 
As we mentioned in \S 5, the $2k$-th moment of zeta should be dominated by  values of $|\zeta(\frac 12+it)|$ of size $(\log T)^k$, and the 
\eqref{6.6}  shows that this set has measure $\ll T(\log T)^{-k^2 +o(1)}$, which yields $M_k(T) \ll T(\log T)^{k^2+\epsilon}$.  

Harper's sharp upper bound for $M_k(T)$ builds on some of these ideas, but deals directly with the moments rather than going through the intermediary of the large deviations in Selberg's theorem \eqref{6.6}.  Instead there is an elaborate decomposition of the sum over primes in \eqref{6.5} into many ranges, and then 
the exponentials of such sums are handled by approximating these by suitable truncations of their Taylor expansion.  Similar ideas were developed independently 
around the same time in \cite{RaSo2} for bounding small moments unconditionally, and the recent paper \cite{HeapSo} develops these ideas in the context of lower bounds.  Thus the proofs of all three aspects of Theorem 6.1 have a unified feel, and the spirit of the proofs may be described as thinking in terms of Euler products but performing computations by replacing Euler products by short Dirichlet series obtained from their Taylor expansions.  These proofs were also influenced  by 
ideas from sieve theory.  For example, in analogy with \eqref{6.5} we may note that $\omega(n)$ (the number of prime factors of $n$) may be bounded above 
by $\sum_{p|n, p\le y} 1  + (\log n)/\log y$ for any $y$, and this could be used to give upper bounds for the mean-value of $d_k(n)$ (which is roughly $k^{\omega(n)}$) 
in short intervals.

The ideas behind obtaining conditional bounds for moments have found diverse applications.  Soundararajan and Young \cite{SoY} used such bounds 
for ``shifted moments'' (see also \cite{Ch1}) to obtain an asymptotic formula (on GRH) for the second moment of quadratic twists of an eigenform $\sum_{|d|\le X} L(\frac 12, f\times\chi_d)^2$.   This is a tantalizing problem, which falls within the purview of the rule of thumb described at the beginning of this section, but 
an unconditional asymptotic has so far been elusive.   A similar problem is to compute the asymptotic for the fourth moment of quadratic Dirichlet $L$-functions 
$\sum_{|d|\le X} L(\frac 12,\chi_d)^4$, and recently Shen \cite{Shen} has extended the method in \cite{SoY} to obtain (on GRH) such an asymptotic.   Analogues of 
these two problems over function fields have been established in \cite{Florea, BFKR}, and since GRH is known in this setting, the corresponding results 
hold unconditionally.  

In a very different direction, Lester and Radziwi{\l \l} \cite{LesRad}  showed on GRH that the Fourier coefficients of half-integer weight Hecke cusp forms 
exhibit a positive proportion of sign changes as we range over fundamental discriminants.   Among the many innovations in their beautiful proof, is an application 
of the ideas discussed above to obtain sharp upper bounds for the second mollified moment of quadratic twists of the Shimura correspondent of the given 
half-integer weight form.   This realization that sharp upper bounds for the second mollified moment suffice has led to another striking result in the 
work of David, Florea, and Lalin \cite{DFL}, who show that a positive proportion of  $L$-functions attached to cubic characters (in the function field setting) have non-zero central value.  Two other recent applications include Zenz \cite{Zenz} to bounding the $L^4$ norm of Hecke eigenforms 
of large weight $k$ for the full modular group, and Shubin \cite{Shubin} to bounding the variance of lattice points on the sphere in random small spherical caps.  
See  \cite{Milin, MilNg, GaoZhao} for further examples.

\section{Extreme values} 

\noindent In \S 2 and \S 3 we discussed the typical size of $|\zeta(\frac 12+it)|$ and central values of $L$-functions, which are governed by Selberg's central limit theorem, and the analogous Keating-Snaith conjectures.  In \S 5 and \S 6 we discussed how the moment problem aims for an understanding of the \emph{large deviations} range of values of $|\zeta(\frac 12+it)|$ (or $|L(\frac 12)|$).  We now discuss the maximal size of $|\zeta(\frac 12+it)|$ (for $0\le t\le T$) and analogous problems in families of $L$-functions.   

As we mentioned in \S 4, our unconditional knowledge is far from the Lindel{\" o}f hypothesis that $|\zeta(\tfrac 12+it)| 
\ll (1+|t|)^{\epsilon}$, and for general $L$-functions already the subconvexity problem poses formidable difficulties.   In 1924 Littlewood established that the Riemann Hypothesis implies the Lindel{\" o}f hypothesis in the quantitative form  
\begin{equation} 
\label{7.1} 
|\zeta(\tfrac 12+it)| \ll \exp\Big( \frac{C \log |t|}{\log \log |t|}\Big) 
\end{equation} 
for some constant $C$.  The estimate \eqref{6.5} yields such a result, upon taking $x =(\log t)^2$ there, and bounding the sum over prime powers trivially.   This strategy was optimized in \cite{ChSo} which showed that one may take any $C> \frac{\log 2}{2}$  in \eqref{7.1}.  Apart from this refinement of the constant $C$, no improvement has been made over Littlewood's estimate.  Corresponding results hold for general $L$-functions, and explicit versions of such bounds (which are useful in computational applications) may be found in \cite{Ch2}. 

Complementing \eqref{7.1}, one may ask for lower bounds on $\max_{T\le t\le 2T} |\zeta(\tfrac 12+it)|$.   Recall that in \S 1 we discussed the extreme values of 
zeta and $L$-functions at the edge of the critical strip, and already there was a gap in our knowledge between the extreme values that may be exhibited 
and the bounds that follow from GRH (see the discussion surrounding \eqref{1.2} and \eqref{1.3}).    This gap becomes much more pronounced on the critical line.   
By using lower bounds for integer moments of $\zeta(\frac 12+it)$, with attention to the uniformity in $k$, Balasubramanian and Ramachandra \cite{BR} 
(optimized in \cite{Balu2}) established that 
\begin{align}
\label{7.2}
\max_{T\le t\le 2T} |\zeta(\tfrac 12+it)| &\ge \max_{k} \Big( \frac{1}{T} \int_T^{2T}  |\zeta(\tfrac 12+it)|^{2k} dt \Big)^{\frac 1{2k}}\nonumber \\
&\gg \max_{k \in {\mathbb N}} \Big( \sum_{n\le T} \frac{d_k(n)^2}{n} \Big)^{\frac 1{2k}} = \exp\Big( (B+o(1)) \frac{\sqrt{\log T}}{\sqrt{\log \log T}}\Big), 
\end{align} 
with $B \approx 0.53$.  With the development of lower bounds for moments in families of $L$-functions (discussed in \S 6), such bounds also became 
available for central $L$-values.  However, a different \emph{resonance method} developed in \cite{So3} has proved to be still more efficient.  The main idea 
in \cite{So3} is to find a Dirichlet polynomial $R(t) = \sum_{n} r(n) n^{-it}$ which ``resonates'' with $\zeta(\frac 12+it)$ and picks out its large values.  This is 
based on computing 
\begin{equation} 
\label{7.3} 
I_1 = \int_T^{2T} |R(t)|^2 dt, \qquad \text{ and  } I_2 = \int_{T}^{2T} \zeta(\tfrac 12+it) |R(t)|^2 dt, 
\end{equation} 
and noting that 
\begin{equation}
\label{7.4} 
\max_{T \le t \le 2T} |\zeta(\tfrac 12+it)| \ge \frac{|I_2|}{I_1}. 
\end{equation} 
If the resonator Dirichlet polynomial $R(t)$ is short, in the sense that $r(n)=0$ unless $n\le T^{1-\epsilon}$, then $I_1$ and $I_2$ in 
\eqref{7.3} may be evaluated asymptotically, and these quantities give two quadratic forms in the unknown coefficients $r(n)$.  The ratio of 
these two quadratic forms is maximized in \cite{So3}, yielding 
\begin{equation} 
\label{7.5} 
\max_{T \le t\le 2T} |\zeta(\tfrac 12+it)| \ge \exp\Big( (1+o(1)) \frac{\sqrt{\log T}}{\sqrt{\log \log T}}\Big). 
\end{equation} 
While this is only a little bit better than \eqref{7.2}, the method also yields lower bounds on the measure of the set on which large values are attained: 
\begin{equation} 
\label{7.6} 
\text{meas}\Big\{ t\in [T, 2T]: \ |\zeta(\tfrac 12+it)| \ge e^V \Big\} \gg \frac{T}{(\log T)^4}  \exp \Big( -10 \frac{V^2}{\log \frac{\log T}{8V^2\log V}}\Big), 
\end{equation} 
uniformly for $3\le V \le \frac 15 \sqrt{\log T/\log \log T}$.  There is some scope to improve such bounds, especially when $V$ is of size $C\log \log T$, 
where one would like to match the upper bound in \eqref{6.6} which would be in keeping with Selberg's theorem (see \cite{HeapSo} for more precise
results when $V\le (2-\epsilon) \log \log T$).  The estimate \eqref{7.6} shows that large values on the scale of \eqref{7.5} occur fairly often (on a 
set of measure $ \ge T^{1-C/\log \log T}$)  suggesting that still larger values might exist.   Furthermore, the resonance method extends readily to families of $L$-functions, and thus we may show (for example) that 
\begin{equation} 
\label{7.7} 
\max_{X\le |d|\le 2X} L(\tfrac 12, \chi_d) \ge \exp\Big( \Big(\frac{1}{\sqrt{5}} + o(1)\Big) \frac{\sqrt{\log X}}{\sqrt{\log \log X}}\Big), 
\end{equation} 
and that, for any Hecke eigenform $f$  
\begin{equation} 
\label{7.8} 
\max_{X \le |d| \le 2X}  L(\tfrac 12, f\times \chi_d) \ge \exp\Big( c\frac{\sqrt{\log X}}{\sqrt{\log \log X}}\Big), 
\end{equation}
for a suitable positive constant $c$.  Indeed the large values in \eqref{7.7} and \eqref{7.8} are attained for more than $X^{1-\epsilon}$ discriminants 
$d$ with $X\le |d|\le 2X$.   By Waldspurger's formula, the large values produced in \eqref{7.8} show that fundamental Fourier coefficients of half-integer weight 
eigencuspforms must get large, and the resonance method has been adapted in \cite{GKS} to show that this holds more generally for half-integer 
weight cusp forms (not necessarily an eigenform).  Another application of this resonance method may be found in the work of Milicevic \cite{Mil} who obtains 
large values of Hecke-Maass cusps forms on arithmetic hyperbolic surfaces. 

Bondarenko and Seip \cite{BonSeip} recently made a breakthrough on this problem, by exhibiting still larger values of $|\zeta(\tfrac 12+it)|$.   The key ingredient is a beautiful result on \emph{GCD sums} or \emph{G{\' a}l sums}:  The problem is to find 
\begin{equation} 
\label{7.9} 
\max_{\substack{ |{\mathcal N}|=N}} \sum_{m, n\in{\mathcal N}} \frac{(m,n)}{\sqrt{mn}},
\end{equation} 
where the maximum is over all $N$ element subsets of the natural numbers.   This elegant combinatorial problem turns out to be closely related to maximizing the ratio of quadratic forms (see \cite{ABS}) 
\begin{equation} 
\label{7.10} 
\max_{|{\mathcal N}| = N} \sup_{\underline{x} \in {\mathbb C}^N \neq 0}  \Big(\sum_{m, n \in {\mathcal N}} x_m \overline{x_n} \frac{(m,n)}{\sqrt{mn}}\Big) 
\Big{/}  \Big({\sum_{n} |x_n|^2}\Big). 
\end{equation} 
Bondarenko and Seip \cite{BonSeip, BonSeip2} established that the maximum in \eqref{7.9} (and also \eqref{7.10}) lies between 
$$
N \exp\Big( (1-\epsilon)\frac{ \sqrt{\log N \log \log \log N}}{\sqrt{\log \log N}}\Big) \text{ and } N  \exp\Big( (7+\epsilon)\frac{ \sqrt{\log N \log \log \log N}}{\sqrt{\log \log N}}\Big),
$$ 
De la Bret{\` e}che and Tenenbaum \cite{dlBT} refined this to show that the maximums in \eqref{7.9} and \eqref{7.10} equal  
\begin{equation} 
\label{7.11} 
N \exp\Big( (2\sqrt{2} + o(1))\frac{ \sqrt{\log N \log \log \log N}}{\sqrt{\log \log N}}\Big). 
\end{equation}

The relevance of the bounds for (related) GCD sums to large values of $|\zeta(\sigma+it)|$ was first appreciated by Aistleitner \cite{Ais} who showed that for fixed $\sigma\in (\tfrac 12, 1)$ and $T$ large one has (for some $c_\sigma >0$)
$$ 
\max_{0< t \le T} |\zeta(\sigma+it)| \ge \exp\Big( \frac{c_\sigma (\log T)^{1-\sigma}}{(\log \log T)^{\sigma}}\Big), 
$$ 
 which improved upon earlier applications of the resonance method (see \cite{Vor2, Hil}) but only matched the results obtained by Montgomery \cite{Mo2} using very different ideas (see also \cite{AMM} for large values on the $1$-line, and \cite{AMMP} for analogous results for Dirichlet $L$-functions).  On the critical line, Bondarenko and Seip \cite{BonSeip} obtained a substantial improvement over the previously known large values of $|\zeta(\frac 12+it)|$ (see \eqref{7.5}) by establishing that 
\begin{equation} 
\label{7.12} 
\max_{0< t\le T} |\zeta(\tfrac 12+it)| \ge \exp\Big( (c+o(1)) \frac{\sqrt{\log T \log \log \log T}}{\sqrt{\log \log T}}\Big), 
\end{equation} 
for a positive constant $c$ (in \cite{BonSeip} $c=1/\sqrt{2}$ is permissible, while \cite{dlBT} allows for the improved $c=\sqrt{2}$).  The key insight is that in the resonance method one can choose ``long resonators'' where $R(t)$ is no longer constrained to be a short Dirichlet polynomial ($r(n) =0$ unless $n\le T^{1-\epsilon}$) 
but instead $R(t)$ is allowed to have $T^{1-\epsilon}$ non-zero coefficients $r(n)$ so long as these are \emph{positive}.  This leads to an optimization problem closely related to the GCD/G{\' a}l sums discussed above, and permits the stronger bound in \eqref{7.12}.   Why is it possible to take such long resonators?  Consider a smooth non-negative function $\Phi$ whose Fourier transform ${\widehat \Phi}$ is also non-negative; for example, we could take $\Phi(t) = e^{-t^2/2}$.  In place of  $I_1$ and $I_2$ in \eqref{7.3} consider the smoothed integrals 
\begin{equation} 
\label{7.13} 
\int_{-\infty}^{\infty} |R(t)|^2 \Phi(t/T) dt, \qquad \text{ and  } \qquad \int_{-\infty}^{\infty} \zeta(\tfrac 12+it) |R(t)|^2 \Phi(t/T) dt. 
\end{equation} 
Replacing $\zeta(\tfrac 12+it)$ with its approximation $\sum_{k \le T} k^{-\frac 12- it}$, the second quantity above is approximately 
$$ 
\sum_{k\le T} \frac{1}{\sqrt{k}} \sum_{m, n} r(m) r(n) \int_{-\infty}^{\infty} \Big(\frac{n}{mk} \Big)^{it} \Phi(t/T) dt 
= T \sum_{k\le T} \frac{1}{\sqrt{k}} \sum_{m, n} r(m) r(n) {\widehat \Phi}(T \log (n/mk)). 
$$ 
Since $m$ and $n$ may be much larger than $T$, we are unable to restrict just to the ``diagonal terms'' $n=mk$, but the crucial point is 
that the positivity of ${\widehat \Phi}$, the resonator coefficients $r(m)$, $r(n)$, and the ``coefficients of $\zeta$'' (namely, the function taking $1$ on 
all positive integers) all allow us to keep any terms that we please on the right side above, and ignore other contributions.  In this way, one can get a 
satisfactory lower bound for the ratio of the quantities in \eqref{7.13}, without needing to evaluate each of these integrals.   The restriction on the number 
of terms allowed in the resonator arises from the fact that $\sum_{k \le T} k^{-\frac 12- it}$ is a poor approximation to $\zeta(\tfrac 12+it)$ if 
$t$ is small.  These small values of $t$ are unavoidable because the condition that  ${\widehat \Phi}$ is non-negative forces
$\Phi(0)$ to be strictly positive.  

 Unlike the resonance method which applies in great generality, there are (at present) limitations on when the Bondarenko--Seip method of using long resonators 
 applies.  In the first place, as we noted above small $t$ must be included, and therefore the maximum in \eqref{7.12} is over 
 $t \in [0,T]$ (this can be refined to the interval $[T^{\beta},T]$ for any $\beta <1$ at the cost of weakening the constant $c$ in \eqref{7.12}), rather than the dyadic intervals $[T,2T]$ seen in \eqref{7.5}.  More significantly, the method requires the positivity of the Dirichlet series coefficients of the $L$-functions in 
 question (analogously to $\zeta$ having coefficients $1$), and also the positivity of the right side of any orthogonality relation or trace formula (analogously to ${\widehat \Phi}$ being non-negative).  Apart from $\zeta(s)$, there is one other example in which the Bondarenko--Seip method has been successfully implemented, and this is the work of de la Breteche and Tenenbaum \cite{dlBT} which produces large values of $|L(\tfrac 12, \chi)|$ as $\chi$ varies over Dirichlet characters $\pmod q$ with $q$ a large prime.   To illustrate the subtleties involved, we note that \cite{dlBT} exhibits large values of $|L(\frac 12,\chi)|$ for \emph{even} characters $\chi$, but the method does not work for \emph{odd} character.  This is because in the even case the orthogonality relation 
 $$ 
 \sum_{\substack {\chi \pmod q \\ \chi \text{even} }} \chi(a) = 
 \begin{cases} 
 \frac{\phi(q)}{2} &\text{ if }  a \equiv \pm 1 \pmod q \\ 
 0 &\text{otherwise} 
 \end{cases}
$$ 
involves only non-negative terms on the right side, whereas this is not the situation for odd characters
$$ 
\sum_{\substack {\chi \pmod q \\ \chi \text{odd} }} \chi(a) = 
 \begin{cases} 
\pm \frac{\phi(q)}{2} &\text{ if }  a \equiv \pm 1 \pmod q \\ 
 0 &\text{otherwise}.
 \end{cases}
$$ 
In particular, the results in \eqref{7.7} and \eqref{7.8} remain the best currently known, and it would be of great interest to see if the Bondarenko--Seip method 
could be extended to more general situations.  

\medskip 

There is a vast gulf between the conditional upper bounds for $|\zeta(\tfrac 12+it)|$ in \eqref{7.1} and the large values exhibited in \eqref{7.12}, and 
it is natural to ask which of these is closer to the truth.  Already in Section 1 we saw a gap (of a factor of $2$) between the extreme values of $L(1,\chi_d)$ 
that may be exhibited (see \eqref{1.2}) and the conditional bounds on these extreme values (see \eqref{1.3}).  There the probabilistic models suggested that the 
extreme values exhibited in \eqref{1.2} represented the truth, and on the critical line too we expect the large values exhibited in \eqref{7.12} to be closer to the 
truth than the bounds in \eqref{7.1}.   For example, if we use Selberg's central limit theorem as a guide and extrapolate, then the measure of $t\in [0,T]$ 
with $|\zeta(\tfrac 12+it)| \ge e^V$ may be expected to be $\ll T\exp(-(1+o(1)) V^2/\log \log T)$ (confer \eqref{6.6}).   If $V = (1+\epsilon) \sqrt{\log T \log \log T}$, 
this measure becomes $\le T^{-\epsilon}$, but one can show that if $|\zeta(\frac 12 +it)|$ attains its maximum for $t\in [0,T]$ at $t=t_0$ 
then in an interval $|t-t_0| \le c/\log T$ its values are at least of size $\frac 12 |\zeta(\tfrac 12+it_0)|$ (see Lemma 2.2 of \cite{FGH}).  This suggests that 
$$ 
\max_{t\in [0,T]} |\zeta(\tfrac 12+it)| \le \exp( (1+o(1)) \sqrt{\log T \log \log T}). 
$$ 
Farmer, Gonek, and Hughes \cite{FGH} have conjectured that even this overestimates the true size of the maximum, and that possibly 
\begin{equation} 
\label{7.14} 
\max_{0\le t\le T} |\zeta(\tfrac 12+it)| = \exp\Big( \Big(\frac{1}{\sqrt{2}}+ o(1)\Big) \sqrt{\log T \log \log T}\Big). 
\end{equation} 
To give one indication of why this might hold, consider \eqref{6.5} which gives (on RH) an upper bound for $\log |\zeta(\tfrac 12+it)|$ in terms 
of essentially a sum over primes going up to $x$, accepting an error term of size $\log T/\log x$.  If we choose $x =\exp(\sqrt{\log T})$ then 
this error term is negligible, and now Re$\sum_{p\le x}  1/p^{\frac 12+it}$ behaves like a Gaussian with mean $0$ and 
variance $\tfrac 12 \sum_{p\le x} 1/p \sim \tfrac 12 \log \log x \sim \tfrac 14 \log \log T$. Extrapolating this Gaussian behavior, we arrive at 
the conjectured behavior in \eqref{7.14}.  The conjecture in \cite{FGH} is based upon a more careful analysis of the hybrid Euler-Hadamard formula 
developed in \cite{GHK}, which decomposes $\log |\zeta(\tfrac 12+it)|$ into terms arising from both primes and zeros in suitable ranges.   Probabilistic 
models for both these terms are analyzed (with the contribution of zeros being modeled using random matrix theory), and the conjecture \eqref{7.14} is consistent with many different ways of splitting into primes and zeros.  Similar conjectures may be formulated in other families of $L$-functions, and for example \cite{FGH} 
conjectures that 
\begin{equation} 
\label{7.15}
\max_{|d|\le X} L(\tfrac 12, \chi_d) = \exp( (1+o(1)) \sqrt{\log X \log \log X}),  
\end{equation} 
which again is a little smaller (by a factor $\sqrt{2}$ in the exponent) than what might be guessed from extrapolating the Keating--Snaith conjectures 
for $\log L(\tfrac 12, \chi_d)$.

As we discussed in \S 4, one motivation for studying the moments of $|\zeta(\tfrac 12+it)|$ is to gain an understanding of its extreme values.  In order to do so, one would need an understanding of how $M_k(T)$ behaves with uniformity in $k$, and a first step might be to examine the asymptotic behavior of the constants $a_k$ and $g_k$ appearing in Conjecture 5.1.  One can show that $\log a_k \sim - k^2 \log (2e^{\gamma} \log k)$, and that $\log g_k \sim k^2 \log (k/4\sqrt{e})$ (see \cite{CoGo}), so that it may seem tempting to speculate that for $T\ge 10$ (say) and uniformly for all $k\ge 2$ one has  (for some positive constant $c$)
$$
T \Big( \frac{c \log T}{k \log k}\Big)^{k^2} \le \int_0^T |\zeta(\tfrac 12+it)|^{2k} dt \le T (\log T)^{k^2}.  
$$ 
But there is a curious paradox, and the upper and lower bounds above are inconsistent!  If the upper bound above holds uniformly, then it follows that 
$$
\max_{0\le t\le T} |\zeta(\tfrac 12+it)| \le \exp( (1+o(1)) \sqrt{\log T \log \log T}).
$$
 Whereas if the lower bound above holds uniformly, then one must 
have 
$$
\max_{0\le t\le T} |\zeta(\tfrac 12+it)| \ge \exp( C \log T/\log \log T)
$$
for some positive constant $C$.  This is an instance where the leading order asymptotic in the moment conjecture does not capture the full story, and one 
should look instead at the recipe in \cite{CFKRS} which (for natural numbers $k$) gives the entire (conjectural) polynomial $P_k$ of degree $k^2$. 
An analysis of this full moment conjecture suggests that the uniform upper bound stated above might hold:  thus, for $T \ge 10$ and natural numbers $k\ge 1$ 
we conjecture that 
\begin{equation} 
\label{7.16} 
\int_0^{T} |\zeta(\tfrac 12+it)|^{2k} dt \le T (\log T)^{k^2}. 
\end{equation}
 In other words, we guess that $\log |\zeta(\tfrac 12+it)|$ is \emph{sub-Gaussian} (when thinking of the frequency of 
 its large values), and this gives a weaker version of the Farmer, Gonek, Hughes conjecture \eqref{7.14}.

 While we have confined our discussion above to large values of $|\zeta(\tfrac 12+it)|$, or equivalently Re$(\log \zeta(\tfrac 12+it))$, similar 
 considerations apply also to Im$(\log \zeta(\frac 12+it))$; see for example \cite{GolGon, CCM, BonSeip3}.

\section{The Fyodorov--Hiary--Keating conjecture} 

\noindent  A fascinating set of problems has emerged recently with the work of Fyodorov and Keating \cite{FK}, and Fyodorov, Hiary, and Keating 
\cite{FHK}, who initiated a study of the distribution of ``local maxima'' of the Riemann zeta function.   More precisely, if $t$ is chosen 
uniformly from $[T,2T]$, what is the distribution of 
$$
\max_{0\le h\le 1} |\zeta(\tfrac 12+it +ih)|?
$$
Although it does not make much of a difference, \cite{FHK} considers the maximum over intervals of length $2\pi$ instead of $1$ since this has a 
natural analogue in random matrix theory.  If a matrix $g$ is chosen randomly from $U(N)$ (with respect to Haar measure), what is the distribution of 
$$ 
\max_{\theta \in [0, 2\pi)} |\text{det} (I - ge^{-i\theta})| ? 
$$ 
In the context of $\zeta(\tfrac 12+it)$, one initial motivation for considering this problem was that it might shed new light on the global 
maximum over the long interval $[0,T]$ (discussed in the previous section).    While the distribution of the local maxima leads to striking new 
and subtle phenomena involving the local correlations of the zeta function, it does not seem to inform the behavior of the global maximum.  

\begin{conjecture} [Fyodorov--Hiary--Keating \cite{FHK}]  For any real number $y$, as $T\to \infty$ one has 
\begin{equation} 
\label{8.1} 
\frac 1T \text{meas} \Big\{ T\le t\le 2T: \ \ \max_{0\le h\le 1} |\zeta(\tfrac 12+it+ih)| \le e^y \frac{\log T}{(\log \log T)^{\frac 34}} \Big\} \to F(y), 
\end{equation} 
where the cumulative distribution function $F$ satisfies $F(y) \to 0$ as $y\to -\infty$, and satisfies $1-F(y) \sim Cye^{-2y}$ as $y\to \infty$ for 
some constant $C>0$.  In particular, 
for any function $g(T)$ tending to infinity with $T$ one has 
\begin{equation} 
\label{8.2} 
 \text{meas} \Big\{ T\le t\le 2T: \ \ \Big| \max_{0\le h\le 1} \log  |\zeta(\tfrac 12+it+ih)| - \log \log T + \frac 34\log \log \log T\Big|\le g(T) \Big\} 
\sim T. 
\end{equation} 
\end{conjecture} 

Let us first explain what is striking and unexpected about this conjecture.  Roughly speaking, in an interval of length $1$ we may think of the 
zeta function as being determined by about $\log T$ values --- this is about the number of zeros we expect to find in such an interval, 
and we may guess that if $|t_1 -t_2| \le 1/\log T$ then $\log |\zeta(\tfrac 12+it_1)|$ and $\log |\zeta(\tfrac 12+it_2)|$ are about the same.   Selberg's theorem 
tells us that the values $\log |\zeta(\tfrac 12+it)|$ are distributed like a normal variable with mean $0$ and variance $\frac 12\log \log T$.  Thus a first guess for the 
distribution of $\max_{0\le h\le 1} \log |\zeta(\tfrac 12+it +ih)|$ might be that it behaves like the maximum of about $\log T$ independently drawn normal random 
variables with mean $0$ and variance $\tfrac 12 \log \log T$.   The maximum of $N$ independent normal variables with mean $0$ and variance $1$ is 
very sharply concentrated around $\sqrt{2\log N} - (\log \sqrt{4\pi \log N})/\sqrt{2\log N}$ (the precise distribution is known as the Gumbel distribution, and has been extensively studied in view of its enormous significance in practical assessments of the risk of rare events).  After scaling by the standard deviation $\sqrt{\frac12 \log \log T}$ in Selberg's theorem, this naive model would indicate that $\max_{0\le h\le 1} \log |\zeta(\tfrac 12+it +ih)|$ should typically be around 
$$ 
\log \log T - \frac{1}{4} \log \log \log T +O(1).
$$
In contrast, Conjecture 8.1 predicts that $\max_{0\le h\le 1} |\zeta(\tfrac 12+it)|$ is usually a bit smaller, of size $(\log T)/(\log \log T)^{\frac 34}$.  
There is also a subtle difference in the decay of $1-F(y)$ in \eqref{8.1}, which is predicted to decay like $ye^{-2y}$, whereas the Gumbel distribution 
would have predicted a decay rate of $e^{-2y}$.

The flaw in the naive heuristic presented above is that nearby values of the zeta function are not independent, but are correlated.  Suppose $t$ is 
randomly chosen from $[T,2T]$ and $0\le h\le 1$, and consider the covariance of $\log |\zeta(\tfrac 12+it)|$ and $\log |\zeta(\tfrac 12+it +ih)|$.  As in our discussion 
of Selberg's theorem in \S 2, we may often approximate these values by corresponding sums over primes Re$\sum_{p\le x} 1/p^{\frac 12+it}$ and Re$\sum_{p\le x} 
1/p^{\frac 12+it +ih}$ with $x$ a suitable small power of $T$.  If $p$ is small in comparison to $e^{1/h}$ then $p^{ih} \approx 1$, and the corresponding terms 
in our prime sums are strongly correlated.   The terms with $p$ much larger than $e^{1/h}$ are largely uncorrelated, since as $p$ varies in such 
large ranges $p^{ih}$ will become equidistributed on the unit circle.   Thus one may see that 
\begin{equation} 
\label{8.3} 
\frac 1T\int_T^{2T} \log |\zeta(\tfrac 12+it)| \log |\zeta(\tfrac 12+it +ih)| dt 
\sim \frac 12 \sum_{p\le x} \frac{\cos(h\log p)}{p} \sim \frac 12 \log \min \big( h^{-1}, \log T\big). 
\end{equation}
This correlation structure of nearby values must be taken into account when trying to predict the behavior of local maxima. 

To gain a rough idea of how to model the local behavior of $\log |\zeta(\frac 12+it)|$, put for each $1 \le k \le \log \log T -1$
\begin{equation} 
\label{8.4} 
{\mathcal P}_k(t) = \text{Re } \sum_{e^{e^{k-1}} \le p \le e^{e^k} } \frac{1}{p^{1/2 +it}}, 
\end{equation} 
so that we may think of $\log |\zeta(\tfrac 12+it)|$ as something like $\sum_{k} {\mathcal P}_k(t)$.  Each ${\mathcal P}_k(t)$ is 
distributed like a normal random variable with mean $0$ and variance $\sim \frac 12 \sum_{e^{e^{k-1}} \le p \le e^{e^k}} 1/p \sim \frac 12$.  
Moreover for different values of $k$, the sums ${\mathcal P}_k(t)$ involve primes in disjoint ranges, and therefore behave independently of 
each other.   Notice further that if $|t_1 -t_2| \le e^{-k}$ then ${\mathcal P}_k(t_1)$ and ${\mathcal P}_k(t_2)$ are more or less the same.  Thus instead of 
modeling $\log |\zeta(\tfrac 12+it)|$ in intervals of length $1$ by about $\log T$ independent samples of a normal random variable, we are led to the following more 
nuanced model.  For each $k$, let $P_k$ denote any one of about $e^k$ independent drawings of a normal random variable with mean $0$ and variance $\frac 12$.  Then $\log |\zeta(\tfrac 12+it)|$ in an interval of length $1$ is modeled by all the possibilities for  $\sum_k P_k$.

The model described above has been analyzed in the probability literature surrounding branching random walks and branching Brownian motion.  Consider  
a particle starting at time $0$ and moving as a standard Brownian motion.  At time $t$, with probability $e^{-t}$ the particle might split into two particles, that 
move according to independent standard Brownian motions starting from that position.  These particles may again split (independently of each other) at a future time, 
giving rise to more daughter particles, and so on.  After time $T$, how is the maximum value of all these particles distributed?   This problem was 
resolved by Bramson who established that the maximum is almost surely $\sqrt{2} (T - \frac 34 \log T )+ O(1)$.  Notice the $\frac 34$ term here, which exactly parallels the $\frac 34$ terms appearing in Conjecture 8.1!

In recent years there has been a lot of progress towards understanding Conjecture 8.1.   In \cite{ABH} Arguin, Belius, and Harper considered $\max_{0 \le h\le 1} \text{Re} \sum_{p\le T} X(p)/p^{\frac 12+ih}$ where the $X(p)$'s are independent random variables  chosen uniformly on the unit circle (a randomized model for 
$\log |\zeta(\tfrac 12+it+ih)|$), and established 
that almost surely this is $\log \log T - (\frac 34+o(1)) \log \log \log T$.  Najnudel \cite{Naj} established that on RH the set of $t \in [T,2T]$ with 
$\max_{0\le h\le 1} |\zeta(\tfrac 12+it +ih)| = (\log T)^{1+o(1)}$ has measure $\sim T$.  Independently this result was also established unconditionally by 
Arguin, Belius, Bourgade, Radziwi{\l \l} and Soundararajan \cite{ABBRS}.   A lovely exposition of Conjecture 8.1 and the results mentioned so far may be 
found in Harper's Bourbaki seminar \cite{Harper4}.  Still more recently, Harper \cite{Harper2} established that if $t$ is not in an exceptional subset of $[T,2T]$ 
with measure $o(T)$, then 
$$ 
\max_{0\le h\le 1} \log |\zeta(\tfrac 12+it+ih)| \le \log \log T -\frac 34 \log \log \log T +O(\log \log \log \log T),
$$ 
so that at least in one direction, the difference between the naive constant $\frac 14$ and the refined prediction $\frac 34$ could be established.  Independently 
Arguin, Bourgade and Radziwi{\l \l} \cite{ABR} established the shaper result that for any $y \ge 1$ 
$$ 
\frac 1T \text{meas} \Big\{ t\in [T,2T]: \ \max_{0\le h\le 1} |\zeta(\tfrac 12+it+ih)| > \frac{e^y \log T}{(\log \log T)^{\frac 34}}\Big\} \le C y e^{-2y}, 
$$ 
for some constant $C$.  This beautiful result establishes part of Conjecture 8.1, and the decay in $y$ above matches (up to constants) 
the conjectured behavior of $1-F(y)$.  There has also been substantial progress toward the analogue of Conjecture 8.1 in random matrix theory; 
see \cite{PaqZeit,  ABB, CMN}.  

Instead of considering the maximum of the zeta function in intervals of length $1$, one may also examine other ``local moments'' 
$\int_0^1 |\zeta(\tfrac 12+it+ih)|^{\beta} dh$.  This was already suggested in \cite{FHK}, who conjectured that a transition in the behavior of these local 
moments occurs at the critical exponent $\beta=2$ --- for $\beta <2$ these local moments are typically of size $(\log T)^{\beta^2/4}$ (the size of the global moment $\frac 1T\int_T^{2T} |\zeta(\tfrac 12+it)|^{\beta} dt$), whereas for $\beta >2$ they are typically of size $(\log T)^{\beta-1}$ corresponding to the largest value of zeta in that interval (about size $\log T$) which might be expected to occur on an interval of length about $1/\log T$.   For work in this direction see \cite{BaiKea2, AOR, Harper2}.  We mention a 
lovely result of Harper \cite{Harper2} for the critical exponent $\beta=2$: 
$$ 
\frac 1T \int_T^{2T} \Big( \frac{1}{\log T} \int_0^1 |\zeta(\tfrac 12+it+ih)|^2 dh \Big)^{\frac 12} dt \ll \frac{1}{(\log \log T)^{\frac 14}}.   
$$ 
A simple application of Cauchy's inequality together with the second moment of $\zeta(\frac 12+it)$ shows that the above quantity is $\ll 1$,   and the 
fact that it is a little bit smaller is a reflection of the correlation structure of nearby values of $\zeta(s)$ that also underlies Conjecture 8.1.

The ideas discussed here are closely connected to what is termed \emph{Gaussian multiplicative chaos}, which was initiated by Kahane \cite{Kahane}, and which has been extensively studied in the probability literature \cite{RhodesVargas}.  In number theory, these ideas are closely related to the study of mean values of random multiplicative functions.  We content ourselves with giving a few pointers to surveys and related work: 
\cite{BaiKea, Harper1, SaksmanWebb, SoundZaman}.

%------
% Insert acknowledgments and information
% regarding funding at the end of the last
% section, i.e., right before the bibliography.
%------

\noindent{\bf Acknowledgements.}  
I am grateful to Brian Conrey, Jon Keating, Emmanuel Kowalski, Vivian Kuperberg, Maksym Radziwi{\l \l}, Matt Tyler, and Max Xu for their careful reading  and many valuable suggestions. This work was partially supported by grants from the National Science Foundation, and a Simons Investigator Award from the Simons Foundarion.

%------
% Insert the bibliography.
%------

\bibliographystyle{emss}
\bibliography{MomRefs}{}

%------ Example for a paper in journal:
% \bibitem{article1}
% A.~Petrunin, Parallel transportation for Alexandrov space with curvature bounded below.
% \emph{Geom. Funct. Anal.} \textbf{8} (1998), no.~1, 123--148.

%------ Example for a book:
% \bibitem{book1}
% W.~P. Ziemer, \emph{Weakly differentiable functions}.
% Grad. Texts in Math. 120,  Springer, New York, 1989.

%------ Example for a paper in a book:
% \bibitem{incollection1}
% J.~S. Milne, Introduction to Shimura varieties.
% In \emph{Harmonic analysis, the trace formula, and Shimura varieties},
% edited by M.~W. Marcellin and E.~Giorgi, pp. 265--378,
% Clay Math. Proc. 4, Amer. Math. Soc., Providence, RI, 2005.

%------ Example for a preprint on arXiv:
% \bibitem{preprint1}
% D.~V. Nguyen, S.~K. Chilappagari, M.~W. Marcellin, and B.~Vasic,
% LDPC codes from latin squares free of small trapping sets,
% 2010, \href{http://arxiv.org/abs/1008.4177}{arXiv:1008.4177}.

%------ Example for a report:
% \bibitem{report1}
% J.~Schöberl, Commuting quasi-interpolation operators.
% Technical report isc-01-10-math, Texas A\&M University, 2001,
% \url{www.isc.tamu.edu/publications-reports/tr/0110.pdf}.

%------ Example for a thesis:
% \bibitem{thesis1}
% E.~Giorgi, \emph{The geometric universe}.
% Ph.D. thesis, University of Maryland, College Park, 2002.

%\end{thebibliography}

\end{document}